\documentclass[preprint,12pt]{elsarticle}

\usepackage{amssymb}
\usepackage{amsmath}
\usepackage{amsthm}

\usepackage[percent]{overpic}
\usepackage{bm,color}
\usepackage{epsfig,verbatim,epstopdf,graphics}
\usepackage{changebar}
\usepackage{multirow}
\usepackage{enumerate}
\usepackage{yhmath}
\usepackage{booktabs} 
\usepackage{cleveref}

\usepackage{tikz}
\usepackage{verbatim}
\usepackage{appendix}
\usetikzlibrary{arrows,backgrounds,decorations.pathreplacing,shapes}
\numberwithin{equation}{section}

\usepackage{subcaption}

\graphicspath{{./}{./figure/}}
\allowdisplaybreaks

\theoremstyle{definition} 

\newtheorem{thm}{Theorem}[section]
\newtheorem{cor}[thm]{Corollary}

\newtheorem{defn}[thm]{Definition}
\newtheorem{rem}[thm]{Remark}

\newtheorem{exa}[thm]{Example}

\newcommand{\norm}[1]{\left\lvert \left\lvert #1 \right\rvert \right\rvert}

\newcommand{\bF}{{\bF F}}

\newcommand{\be}{\begin{eqnarray}}
\newcommand{\ee}{\end{eqnarray}}
\newcommand{\beno}{\begin{eqnarray*}}
\newcommand{\eeno}{\end{eqnarray*}}

\makeatletter

\newcommand{\Rmnum}[1]{\expandafter\@slowromancap\romannumeral #1@}
\makeatother

\journal{Journal of Computational Physics}

\begin{document}

\begin{frontmatter}



\title{An Adaptive-rank Approach with Greedy Sampling 
for Multi-scale BGK Equations 
}


\author[inst1]{William A. Sands}
\ead{wsands@udel.edu}

\author[inst1]{Jing-Mei Qiu}
\ead{jingqiu@udel.edu}

\author[inst1]{Daniel Hayes}
\ead{dphayes@udel.edu}

\author[inst1]{Nanyi Zheng\corref{cor1}}
\ead{nyzheng@udel.edu}

\cortext[cor1]{Corresponding author}

\affiliation[inst1]{
  organization={Department of Mathematical Sciences,  University of Delaware},
  city={Newark},
  postcode={19716},
  state={Delaware},
  country={USA}
}
\begin{abstract}
In this paper, we propose a novel adaptive-rank method for simulating multi-scale BGK equations which is based on a greedy sampling technique. The method adaptively identifies important rows and columns of the solution matrix, and reduces computational complexity by updating only the selected rows and columns through a local solver. Once updated, the solution at selected rows and columns is used together with an adaptive cross approximation to reconstruct the entire solution matrix. The approach extends the semi-Lagrangian adaptive-rank approach, introduced in our previous work for the Vlasov-Poisson system \cite{zheng2025semi} in several ways. Unlike the step-and-truncate low-rank integrators discussed in \cite{einkemmer2025review}, the greedy sampling technique considered in this paper avoids the need for explicit low-rank decompositions of nonlinear terms, such as the local Maxwellian in the BGK collision operator. We enforce mass, momentum, and energy conservation by introducing a new locally macroscopic conservative correction, which implicitly couples the kinetic solution to the solution of the corresponding macroscopic system. Using asymptotic analysis, we show that the macroscopic correction preserves the asymptotic properties intrinsic to the full-grid scheme and that the proposed method in the low-rank setting possesses a conditionally asymptotic-preserving property. Another unique advantage of our approach is the use of a local semi-Lagrangian solver, which permits large time steps compared to Eulerian schemes. This flexibility is retained in the macroscopic solver by employing high-order stiffly-accurate diagonally implicit Runge-Kutta methods. The resulting nonlinear macroscopic systems are solved efficiently using a Jacobian-free Newton-Krylov method, which eliminates the need for preconditioning at modest CFL numbers. Each iteration of the nonlinear macroscopic solver provides a self-consistent correction to a provisional low-rank kinetic solution, which is then used as a dynamic closure for the macroscopic system. Numerical results demonstrate the efficacy of the proposed method in capturing shocks and discontinuous solution structures. We also highlight its performance in a challenging mixed-regime problem, where the Knudsen number spans multiple orders of magnitude.
\end{abstract}




\end{frontmatter}


\section{Introduction}
\label{sec:intro}

The Bhatnagar–Gross–Krook (BGK) equation \cite{bgk1954model} is a popular model used to study the dynamics of rarefied gases when the mean-free path is large compared to characteristic length scales, causing fluid models such as the Euler and Navier Stokes equations to break down. It is a nonlinear relaxation model that approximates the collision dynamics of the well-known Boltzmann collision operator which requires the evaluation of a five-dimensional integral that is prohibitively expensive due to its intrinsic nonlinearity and high-dimensionality. In contrast, the BGK operator is a simpler nonlinear algebraic operator; yet it retains key properties of the Boltzmann operator, such as entropy dissipation and the preservation of collision invariants (e.g., mass, momentum, and energy) and various multiscale phenomena. As a result, the BGK equation serves as a valuable prototype for studying more complex collision operators.

Low-rank tensor approximations in the context of kinetic simulations have attracted significant interest in recent years \cite{einkemmer2025review}. Such methods compress the time-dependent numerical solution using low-rank decomposition techniques developed for matrices and tensors, and are a promising tool to address the challenges associated with the curse of dimensionality. The primary research directions in this area can be interpreted as either dynamic low-rank (DLR) or step-and-truncate (SAT) methods. DLR methods \cite{koch2007dynamical, einkemmer2018low,einkemmer2019low,einkemmer2021efficient,ceruti2022unconventional,baumann2024stable} project the continuous equation onto a low-rank manifold using the low-rank bases from the solution, and then discretize the resulting equations to obtain evolution equations for the bases. In contrast, SAT methods \cite{kormann2015semi, dektor2021rank, GuoVlasovFlowMap2022, guo2024local, sands2025high} first discretize the kinetic equation and dynamically evolve a low-rank representation of the solution using tensor linear algebra. One of the challenges associated with the BGK collision operator in the context of low-rank methods is the need to evaluate a local Maxwellian, which generally does not admit a low-rank decomposition. To avoid these challenges, some of the existing low-rank methods for the BGK equation in the DLR \cite{einkemmer2021efficient,baumann2024stable} framework use a multiplicative splitting that is guided by the structure of the solution obtained from Chapman-Enskog analysis. Interpolatory low-rank methods, also known as CUR or cross approximation \cite{shi2024distributed}, or an oblique projection in the DLR framework \cite{dektor2025nonlinear,dektor2024interpolatoryBGK}, employ a sampling strategy based on the discrete empirical interpolation method (DEIM) that allows one to avoid forming explicit low-rank decompositions. This is particularly advantageous for the evaluation of the local Maxwellian in the BGK equation, where such methods were recently used to investigate compression of the entire six-dimensional distribution function \cite{dektor2024interpolatoryBGK}. In a recent work \cite{zheng2025semi}, the adaptive cross approximation of matrices \cite{shi2024distributed} was used to construct a non-split semi-Lagrangian adaptive-rank (SLAR) method that was applied to the nonlinear Vlasov-Poisson system \cite{zheng2025semi}. 
Compared with existing low-rank approaches, the sampling-based adaptive-rank methods uses information from a residual to adaptively identify critical rows and columns that need to be updated. Once the rows and columns are identified, the low-rank approximation is then updated using low-cost recursive rank-one corrections. These rows and columns are a subset of the entire matrices, with a computational complexity similar to that of a low-rank solver. The method is also demonstrated to be highly flexible in the development of numerical schemes, can achieve high-order convergence in spatial and temporal discretizations, and can achieve the mass conservation property via a locally macroscopic conservation (LoMaC) type correction, similar to those in \cite{guo2024local, guo2024conservative}.    

The goal of this work is to develop a new adaptive-rank approach for the multi-scale BGK equation that employs a greedy sampling strategy. The local solver used in this work is built upon a SL finite difference method \cite{li2023high}, which achieves high-order accuracy in both temporal and spatial discretizations and allows for rather large time steps with numerical stability. In \cite{ding2023accuracy}, such a local solver was shown to be asymptotic preserving (AP) and asymptotically accurate (AA) for the BGK equation in the semi-discrete sense. This work is built upon such a local solver, but with the greedy sampling technique and the SLAR method from our recent work \cite{zheng2025semi}, which was shown to be flexible in its implementation, robust in its numerical performance, and efficient with a low-rank complexity. A key advantage of the adaptive-rank approach considered in this work is that explicit low-rank decomposition of local Maxwellians are no longer necessary, which greatly simplifies the construction of low-rank solvers for the BGK equation. This work also introduces a new LoMaC correction which can simultaneously achieve mass, momentum, and energy conservation. The original version of LoMaC proposed in \cite{guo2024local} considered explicit schemes for the macroscopic system. To eliminate the CFL restriction from the macroscopic system, we consider implicit discretizations that couple the macroscopic and microscopic systems in order to define a dynamic and self-consistent correction of the low-rank kinetic solution. The approach we use is inspired by HOLO methods \cite{taitano2014moment,taitano2015charge,park2019multigroup,hammer2019multi} which have shown great success in the simulation of highly nonlinear kinetic problems. At its core, the HOLO method uses a macroscopic (moment) formulation with a dynamic closure to accelerate the convergence of a nonlinear scheme for the kinetic equation. The nonlinear schemes for the kinetic equation are often linearized using Picard iteration, which requires one to solve the kinetic equation multiple times. In contrast, the approach we present first computes a high quality provisional low-rank solution with a non-conservative SL method. Then, using a nonlinear iterative method, we simultaneously correct the provisional solution in a way that directly enforces the macroscopic conservation laws. In order to retain the efficiency of the SL method, we discretize the macroscopic equations using high-order stiffly-accurate (SA) diagonally-implicit Runge-Kutta (DIRK) methods. The resulting nonlinear systems are solved using a Jacobian-free Newton-Krylov (JFNK) method, which, for modest CFL, does not require preconditioning and uses the high-quality non-conservative provisional solution as an initial guess. We also would like to remark that the proposed LoMaC approach is quite general and does not require a low-rank decomposition of nonlinear terms.

The organization of the paper is structured as follows. First in \Cref{sec:models}, we introduce the kinetic model and its associated macroscopic system. The bulk of this paper is concentrated in \Cref{sec:discretizations} and subsections therein. We provide a brief overview of the high-order SL finite difference method for the BGK equation \cite{li2023high} in \Cref{subsec:deterministic_SLFD}, which is the full rank version of the proposed adaptive-rank method described in \Cref{subsec:SLAR for the BGK equation}. In \Cref{subsec:moment preservation}, we provide details regarding the preservation of the macroscopic conservation laws, and discuss the techniques used to solve the coupled nonlinear macroscopic system. \Cref{subsec:high-order DIRK} describes extensions of the aforementioned methods to high-order temporal accuracy using DIRK. We also conduct an asymptotic analysis of the method near the fluid limit and prove that the scheme possesses a \textit{conditionally} AP property. Numerical results are provided in \Cref{sec:numerical_tests} to confirm the accuracy and robustness of the proposed methodologies. \Cref{sec:conclusion} contains the conclusion, which provides an overview of the results and discusses extensions for future work.

\section{The Multi-scale Kinetic Equation and Macroscopic System}
\label{sec:models}

In this section, we provide details of the models that will be investigated in this work as well as the conventions used regarding notation. We first introduce the kinetic model in \Cref{subsec:BGK model} and discuss its corresponding macroscopic system in \Cref{sec:fluid model}.

\subsection{The BGK Equation}
\label{subsec:BGK model}
We consider the BGK relaxation model \cite{bgk1954model} of the Boltzmann equation, which, for neutral particles, in the absence of external forces, leads to the kinetic model
\begin{equation}
    f_t + vf_x = \frac{\mathcal{M}_{\mathbf{U}(f)} - f}{\epsilon}, \quad \left(x, v\right) \in \Omega_x \times \Omega_v. \label{eq:BGK}
\end{equation}
Here, the distribution function $f = f(x,v,t)$ represents the probability of finding a particle at position $x \in \Omega_x$, with velocity $v \in \Omega_v$, at time $t$. The parameter $\epsilon > 0$ is the Knudsen number whose inverse characterizes the relaxation rate of the system to its local thermodynamic equilibrium. As a consequence of the Boltzmann $\mathcal{H}$-theorem, this equilibrium is defined to be the local Maxwellian distribution
\begin{equation}
    \label{eq:local maxwellian}
    \mathcal{M}_{\mathbf{U}(f)} := \frac{\rho}{\sqrt{2\pi T}}\exp\left(-\frac{(v - u)^2}{2T}\right),
\end{equation}
where we have defined the raw moments
\renewcommand{\arraystretch}{1.2}
\begin{align}
    \label{eq:raw moments}
    \rho = \left\langle f \right\rangle_v, \quad \rho u = \left\langle vf \right\rangle_v, \quad E = \left\langle \frac{1}{2}v^{2} f \right\rangle_v, \quad \mathbf{U}(f):=
    \begin{pmatrix}
    \rho \\
    \rho u \\
    E
    \end{pmatrix},
\end{align}
as well as the central moment
\begin{align}
    \label{eq:central moments}
    \frac{\rho T}{2} = \left\langle \frac{1}{2}(v-u)^{2} f\right\rangle_v.
\end{align}
To simplify the notation, we write
\begin{equation}
    \label{eq:velocity integration}
    \left\langle g \right\rangle_v := \int_{\Omega_{v}} g(x,v,t) \,dv,
\end{equation}
to denote integration of a function $g = g(x,v,t)$ in velocity. It is useful to note that the total energy density $E$ can be written in terms of the bulk kinetic energy and internal energy, namely
\begin{equation}
    \label{eq:relation between energy and temperature}
    E = \frac{1}{2}\left( \rho u^{2} + \rho T \right).
\end{equation}
Finally, an important property of the BGK collision operator is that it possesses several invariants. In particular, it preserves mass, momentum, and energy in the sense that
\begin{equation}
    \label{eq:BGK invariants}
    \left\langle \left( \frac{\mathcal{M}_{\mathbf{U}(f)} - f}{\epsilon} \right) \Phi \right\rangle_v = 0, \quad \Phi(v) = \left(1, v, \frac{1}{2}v^{2}\right)^{\top}. 
\end{equation}

\subsection{Macroscopic System}
\label{sec:fluid model}

In addition to the kinetic equation \eqref{eq:BGK}, this work also utilizes the conserved macroscopic system, which is obtained in a self-consistent manner through moments of the distribution function $f$. This system of moment equations can be obtained from \eqref{eq:BGK} by multiplying both sides with a function $\Phi(v) = \left(1, v, \frac{1}{2}v^{2}\right)^{\top}$ and then applying the operator $\left\langle \cdot \right\rangle_{v}$ to both sides. By \eqref{eq:BGK invariants}, the functions in $\Phi(v)$ are known to be collision invariants, so we obtain a system for the macroscopic variables of the form
\begin{equation}
    \label{eq:Fluid system form}
    \mathbf{U}_{t} + \mathbf{F}(\mathbf{U})_{x} = 0,
\end{equation}
where the conserved quantities $\mathbf{U} = \mathbf{U}(f)$ are defined according to \eqref{eq:raw moments}, and the macroscopic flux $\mathbf{F}(\mathbf{U})$ can be written as
\renewcommand{\arraystretch}{1.5}
\begin{equation}
    \label{eq:macro flux F}
  \mathbf{F}(\mathbf{U}):=
    \begin{pmatrix}
    \left\langle vf \right\rangle_v \\
    \left\langle v^{2}f \right\rangle_v \\
    \left\langle \frac{1}{2} v^{3}f \right\rangle_v
    \end{pmatrix}.
\end{equation}
We remark that the macroscopic flux \eqref{eq:macro flux F} can be equivalently written in terms of the conserved variables $\mathbf{U}$, with the flux of the last equation requiring a closure. The structure of the system \eqref{eq:Fluid system form} characterizes the fluid limit of the kinetic equation in the limit $\epsilon \rightarrow 0$, which should be preserved at the discrete level. In particular, as $\epsilon \rightarrow 0$, the distribution can be well-approximated by a local Maxwellian. A standard perturbation argument shows that $f = \mathcal{M}_{\mathbf{U}(f)}$ to leading order in $\epsilon$, which results in the compressible Euler system, if it is used as the closure. Corrections to $\mathcal{O}(\epsilon)$ and beyond can be achieved using the classical Chapman-Enskog expansion \cite{chapman1990mathematical}, leading to the Navier-Stokes and Burnett equations. In this work, we take a slightly different approach in which an approximate solution of the kinetic equation defines the closure of the moment system. This approximate solution shall be obtained using the adaptive-rank method described in the following section, which is developed using a greedy sampling technique. A critical aspect of our approach is that the solution of the moment system is treated in a manner that is self-consistent with the kinetic equation. This, in turn, allows us to simultaneously satisfy key conservation laws and preserve the structure of the fluid limit. We provide additional details of this technique in \Cref{subsec:moment preservation}.

\section{Adaptive-rank Method and Preservation of Conservation Laws}
\label{sec:discretizations}

We begin with a brief review of high-order deterministic (full-rank) SL-FD methods for the BGK equation in \Cref{subsec:deterministic_SLFD}. Then, in \Cref{subsec:SLAR for the BGK equation} and \Cref{subsec:moment preservation}, we describe the proposed sampling-based adaptive-rank method and the iterative procedure used to enforce macroscopic conservation laws. These methods are all developed under a first-order backward Euler time discretization to highlight the key algorithmic ideas. In \Cref{subsec:high-order DIRK}, we extend the approach to high-order accuracy in time using stiffly-accurate DIRK methods. Finally, in \Cref{subsec:asymptotic analysis}, we perform an asymptotic analysis of the proposed scheme.

\subsection{SL-FD Method for the BGK Equation}
\label{subsec:deterministic_SLFD}

In this paper, we discretize the phase space $\Omega_{x} \times \Omega_{v} = [a_{x}, b_{x}] \times [a_{v}, b_{v}]$ with a uniform Cartesian mesh consisting of $N_x \times N_v$ cells whose boundaries are
\begin{align*}
    &a_{x} = x_{\frac{1}{2}} < x_{\frac{3}{2}} < \ldots < x_{N_x+\frac{1}{2}} = b_{x}, \\
    &a_{v} = v_{\frac{1}{2}} < v_{\frac{3}{2}} < \ldots < v_{N_v+\frac{1}{2}} = b_{v}.
\end{align*}
For brevity, we denote the cells, as well as their respective centers and resolutions, as
\begin{align*}
    &I^x_i = [x_{i-\frac{1}{2}},x_{i+\frac{1}{2}}], \quad x_i = \frac{1}{2}(x_{i-\frac{1}{2}}+x_{i+\frac{1}{2}}), \quad \Delta x = x_{i+\frac{1}{2}} - x_{i-\frac{1}{2}}, \quad \text{for } i = 1, \cdots, N_{x}, \\
    &I^v_j = [v_{j-\frac{1}{2}},v_{j+\frac{1}{2}}], \quad v_j = \frac{1}{2}(v_{j-\frac{1}{2}}+v_{j+\frac{1}{2}}), \quad \Delta v = v_{j+\frac{1}{2}} - v_{j-\frac{1}{2}}, \quad \text{for } j = 1, \cdots, N_{v}.
\end{align*}
It is also convenient to write $I_{i,j} = I^{x}_{i}\times I^{v}_{j}$ for all $i,j$ when referring to cells in the two-dimensional grid. Lastly, we can organize the cell centers into vectors, namely
\begin{equation*}
    \mathbf{x} := [x_1, x_2,\ldots, x_{N_x}]^\top, \quad \mathbf{v} := [v_1, v_2,\ldots, v_{N_v}]^\top.
\end{equation*}

Next, using the material derivative, we cast the BGK equation \eqref{eq:BGK} in terms of its characteristics, which gives
\begin{equation}
    \frac{df}{dt} = f_t + vf_x= \frac{\mathcal{M}_{\mathbf{U}(f)} - f}{\epsilon}.
\end{equation}
A first-order backward Euler (BE) time discretization along these characteristics produces
\begin{equation}\label{eq:BE_form1}
    \begin{split}
        f^{n+1}_{i,j} &= f(x_i-v_j\Delta t,v_j,t^n) + \Delta t\frac{(\mathcal{M}_{\mathbf{U}(f^{n+1})})_{i,j} - f^{n+1}_{i,j}}{\epsilon},\\
        &=: \widetilde{f}^{n+1}_{i,j} + \Delta t\frac{\left(\mathcal{M}_{\mathbf{U}(f^{n+1})}\right)_{i,j} - f^{n+1}_{i,j}}{\epsilon}.
    \end{split}
\end{equation}
Next, if we multiply both sides of the update \eqref{eq:BE_form1} by $\Phi(v) = \left(1, v, \frac{1}{2}v^{2}\right)^{\top}$, and apply the operator $\left\langle \cdot \right\rangle_{v}$, with the aid of property \eqref{eq:BGK invariants}, it can be shown that
\begin{equation}
    \label{eq:equivalence of Maxwellians under collision invariance}
    \mathcal{M}_{\mathbf{U}(f^{n+1})} = \mathcal{M}_{\mathbf{U}(\widetilde{f}^{n+1})}.
\end{equation}
Hence, the BE time discretization \eqref{eq:BE_form1} can be rewritten in an explicit form:
\begin{equation}\label{eq:BE_final}
    f^{n+1}_{i,j} = \frac{\epsilon \widetilde{f}^{n+1}_{i,j} + \Delta t \left(\mathcal{M}_{\mathbf{U}(\widetilde{f}^{n+1})}\right)_{i,j}}{\epsilon + \Delta t}.
\end{equation}
To perform the spatial reconstruction of $\widetilde{f}^{n+1}_{i,j}$ from the feet of characteristics as in \Cref{eq:BE_form1}, we use the fifth-order finite difference scheme combined with WENO reconstructions \cite{li2023high}. This BE discretization of the BGK equation was previously shown to be AP and was extended to high order using DIRK methods \cite{ding2023accuracy} with the AA property. In \Cref{subsec:asymptotic analysis}, we extend these ideas to account for the truncation process used in low-rank approximation. As outlined above, we begin by presenting the method under a first-order backward Euler discretization, which serves as the foundation for the adaptive-rank scheme. The extension to high-order accuracy using DIRK methods is discussed separately in \Cref{subsec:high-order DIRK}. 

\subsection{An Adaptive-rank Method for the BGK Equation with Greedy Sampling}
\label{subsec:SLAR for the BGK equation}

To elaborate the idea of the adaptive-rank approach with greedy sampling, we adopt the following notation for the solution matrices \( F^{n} = (f_{i,j}^{n})_{N_x\times N_v} \) and the intermediate solution \( \widetilde{F}^{n+1} = (\widetilde{f}_{i,j}^{n+1})_{N_x\times N_v} \) from the SL update. The element-wise update in \eqref{eq:BE_final} can be rewritten in matrix form as
\begin{equation}\label{eq:BE_final_F}
    F^{n+1} := \textbf{L}(\widetilde{F}^{n+1}) = \textbf{L} \circ \textbf{S} (F^n),
\end{equation}
where \( \textbf{S} \) represents the SL method with reconstruction \cite{li2023high}, and \( \textbf{L} \) denotes the nonlinear update operator defined in \eqref{eq:BE_final}. We next describe the adaptive-rank compression strategy used to efficiently represent each operator in this composition. 

\subsubsection*{ACA+SVD Strategy}

The core compression routine used in our adaptive-rank framework is a two-stage strategy that combines the adaptive cross approximation (ACA) with an SVD truncation. This ACA+SVD strategy, originally proposed as the central component of the SLAR method in our previous work~\cite{zheng2025semi}, enables efficient construction of low-rank approximations using only a subset of matrix entries. The ACA is particularly well-suited for problems where operators, solution updates, or nonlinear terms can be evaluated locally without assembling the full matrix. The subsequent SVD truncation further enhances the robustness of the method by reducing redundancy and improving numerical stability~\cite{zheng2025semi}. We briefly recall the ACA+SVD strategy below, beginning with a visual illustration of the CUR decomposition.

\begin{figure}[ht]
    \centering
    \includegraphics[width=0.65\linewidth]{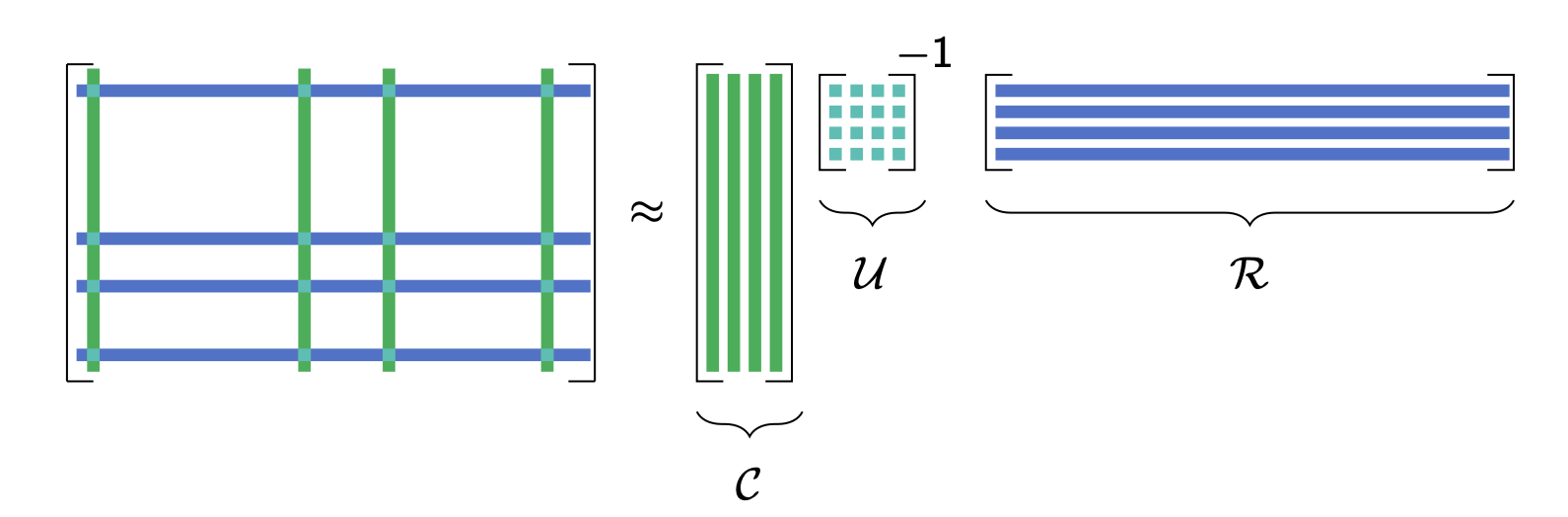}
    \caption{Visual representation of the CUR decomposition constructed by ACA.}
    \label{fig:CUR_Picture}
\end{figure}

Given a matrix \( A \in \mathbb{C}^{N_x \times N_v} \), the ACA algorithm provides an efficient way to construct a CUR decomposition (see \Cref{fig:CUR_Picture}) of the form
\begin{equation}\label{eq:CUR_decom}
A_k = A(:, \mathcal{J})\, A(\mathcal{I}, \mathcal{J})^{-1}\, A(\mathcal{I}, :) =: \mathcal{C}\mathcal{U}\mathcal{R},
\end{equation}
where \( \mathcal{I} = \{i_1, i_2, \ldots, i_k\} \) and \( \mathcal{J} = \{j_1, j_2, \ldots, j_k\} \) are the selected row and column index sets; \( \mathcal{C} \) and \( \mathcal{R} \) denote the corresponding selected columns and rows of \( A \), and \( \mathcal{U} \) is the inverse of the intersection matrix \( A(\mathcal{I}, \mathcal{J}) \). While choosing the optimal index set in a CUR decomposition is NP-hard~\cite{civril2007finding}, ACA builds up the decomposition using a greedy pivoting strategy driven by residual error.

The recursive ACA algorithm~\cite{shi2024distributed,zheng2025semi} constructs the approximation iteratively by updating a rank-\(k\) representation via a rank-one correction to the previous rank-\(k-1\) approximation. Each iteration—referred to as the \(k\)-th step—proceeds as follows:

\paragraph{Phase I: Pivot selection.}
Let \( A_{k-1} \) denote the current approximation (with \( A_0 = (0)_{N_x\times N_v} \)) and \( R_{k-1} = A - A_{k-1} \) the corresponding residual. At iteration \(k\), a candidate set of index pairs \( \mathcal{P} = \{(i_\ell, j_\ell)\}_{\ell=1}^{12} \) is randomly sampled from the unselected rows and columns. The pivot location is initialized as
\(
(i^*, j^*) = \arg\max_{(i,j) \in \mathcal{P}} |R_{k-1}(i,j)|,
\)
followed by greedy refinement steps:
\[
i_k = \arg\max_i |R_{k-1}(i,j^*)|, \quad
j_k = \arg\max_j |R_{k-1}(i_k,j)|.
\]
The pair \((i_k, j_k)\) is then appended to the index sets \( \mathcal{I} \) and \( \mathcal{J} \).

\paragraph{Phase II: Rank-one update.}
A rank-one correction is applied:
\begin{equation}
\label{eq:recursive_update}
A_k = A_{k-1} + \frac{1}{R_{k-1}(i_k, j_k)}\, R_{k-1}(:, j_k)\, R_{k-1}(i_k, :).
\end{equation}
This update is equivalent to the standard CUR decomposition \eqref{eq:CUR_decom}, which ensures interpolation at all selected rows and columns in $\mathcal{I}$ and $\mathcal{J}$: the updated approximation satisfies \( A_k(i_\ell,:) = A(i_\ell,:), \quad A_k(:,j_\ell) = A(:,j_\ell) \quad \text{for } \ell = 1, \dots, k \).  This process continues until the Frobenius norm of the most recent rank-one correction falls below a prescribed tolerance \( \varepsilon_c \), or until a specified maximum rank is reached.

The resulting ACA approximation admits a structured decomposition:
\[
A_k = \mathcal{E}_{\mathcal{J}}\, \mathcal{D}\, \mathcal{E}_{\mathcal{I}},
\]
where \( \mathcal{E}_{\mathcal{J}} = [\mathbf{c}_1, \dots, \mathbf{c}_k] \) with \( \mathbf{c}_\ell = R_{\ell-1}(:, j_\ell) \); \( \mathcal{E}_{\mathcal{I}} = [\mathbf{r}_1, \dots, \mathbf{r}_k]^T \) with \( \mathbf{r}_\ell = R_{\ell-1}(i_\ell, :)^T \); and \( \mathcal{D} = \mathrm{diag}(r_1^{-1}, \dots, r_k^{-1}) \) with \( r_\ell = R_{\ell-1}(i_\ell, j_\ell) \).

\paragraph{SVD truncation}
To improve numerical stability and eliminate redundant modes, the CUR approximation is further compressed using singular value decomposition (SVD). Specifically, we first perform QR decompositions of $\mathcal{E}_{\mathcal{J}}$ and $\mathcal{E}_{\mathcal{I}}^T$, followed by an SVD on the reduced core:
\[
\mathcal{E}_{\mathcal{J}} = Q_1 R_1, \quad \mathcal{E}_{\mathcal{I}}^T = Q_2 R_2, \quad R_1 \mathcal{D} R_2^{\mathrm{T}} = \widetilde{U} \widetilde{\Sigma} \widetilde{V}^{\mathrm{T}}.
\]
We then truncate the decomposition at rank $r_s$, defined by the smallest integer such that $\widetilde{\Sigma}(r_s+1, r_s+1) < \varepsilon_s$, and construct the final approximation:
\[
A \approx \left(Q_1 \widetilde{U}(:,1\!:\!r_s)\right) \widetilde{\Sigma}(1\!:\!r_s,1\!:\!r_s) \left(Q_2 \widetilde{V}(:,1\!:\!r_s)\right)^{\mathrm{T}}=: U\Sigma V^{\mathrm{T}},
\]
which defines the final output format of the ACA+SVD strategy.

We now apply this compression strategy to the operator that defines the BGK update. Using the formulation \eqref{eq:BE_final_F}, we define the corresponding low-rank approximation as
\begin{equation}\label{eq:BE_final_F2}
    \mathcal{F}^{n+1,\star} = \mathcal{L} \circ \mathcal{S} (\mathcal{F}^n),
\end{equation}
where \( \mathcal{S} := \mathcal{T}_{\varepsilon_s} \circ \mathcal{S}_{\varepsilon_c} \) and \( \mathcal{L} := \mathcal{T}_{\varepsilon_s} \circ \mathcal{L}_{\varepsilon_c} \) represent two ACA+SVD stages used to approximate the SL update operator \( \mathbf{S} \) and the nonlinear collision operator \( \mathbf{L} \), respectively. The superscript \( \star \) in \( \mathcal{F}^{n+1,\star} \) indicates that this intermediate result does not yet satisfy local conservation laws. The LoMaC correction procedure, introduced later in \Cref{subsec:moment preservation}, will adjust \( \mathcal{F}^{n+1,\star} \) to enforce conservation of mass, momentum, and energy, yielding the final conservative solution \( \mathcal{F}^{n+1} \).

\begin{figure}[h!]
    \begin{center}
        \begin{overpic}[scale=0.66]{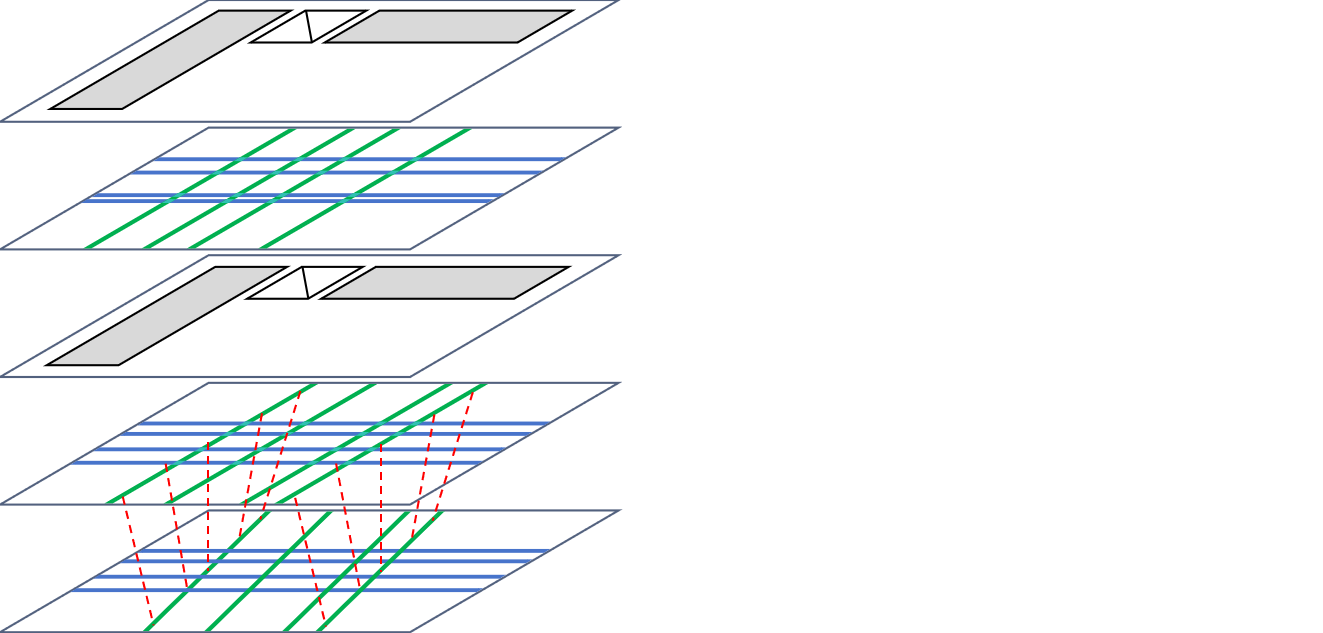}
            \put(48,2) {\small$\mathcal{F}^{n}$}
            
            \put(52,5){
                \begin{tikzpicture}[remember picture, overlay]
                \draw[->, thick, bend right=45, black] (1,0) to (1,1);
                \end{tikzpicture}
            }
            \put(62,7.8) {\footnotesize$\mathcal{S}_{\varepsilon_c}$}
            
            \put(48,12.8) {\small$\widetilde{\mathcal{F}}^{n+1}_{\{\#,k_1\}}$}

            \put(52,15.8){
                \begin{tikzpicture}[remember picture, overlay]
                \draw[->, thick, bend right=45, black] (1,0) to (1,1);
                \end{tikzpicture}
            }
            \put(62,18) {\footnotesize$\mathcal{T}_{\varepsilon_s}$}
            
            \put(48,23.6) {\small$\widetilde{\mathcal{F}}^{n+1}_{\{\#,k_1\!\to r_1\}}$}

            \put(52,26.6){
                \begin{tikzpicture}[remember picture, overlay]
                \draw[->, thick, bend right=45, black] (1,0) to (1,1);
                \end{tikzpicture}
            }
            \put(62,29.4) {\footnotesize$\mathcal{L}_{\varepsilon_c}$}
            
            \put(48,34.4) {\small$\mathcal{F}^{n+1,\star}_{\{\#,k_2\}}$}

            \put(52,37.4){
                \begin{tikzpicture}[remember picture, overlay]
                \draw[->, thick, bend right=45, black] (1,0) to (1,1);
                \end{tikzpicture}
            }
            \put(62,40.2) {\footnotesize$\mathcal{T}_{\varepsilon_s}$}
            
            \put(48,45.2) {\small$\mathcal{F}^{n+1,\star}_{\{\#,k_2\!\to r_2\}}$}

            \put(62,45.2) {\vector(1,0){20}}
            \put(68.5,45.5) {\footnotesize LoMaC}
            \put(85,45) {\small$\mathcal{F}^{n+1}$}

            \put(65,8){
                \begin{tikzpicture}[remember picture, overlay]
                    \draw [decorate, decoration={brace, amplitude=6pt, mirror}, thick] (0,0) -- (0,1.6);
                \end{tikzpicture}
            }
            \put(69,12.8) {\footnotesize$\mathcal{T}_{\varepsilon_s}\circ\mathcal{S}_{\varepsilon_c}=:\mathcal{S}$}

            \put(65,30){
                \begin{tikzpicture}[remember picture, overlay]
                    \draw [decorate, decoration={brace, amplitude=6pt, mirror}, thick] (0,0) -- (0,1.6);
                \end{tikzpicture}
            }
            \put(69,34.8) {\footnotesize$\mathcal{T}_{\varepsilon_s}\circ\mathcal{L}_{\varepsilon_c}=:\mathcal{L}$}
        \end{overpic}
    \end{center}
\caption{Adaptive-rank method for the BGK equation.
}
\label{fig:SamplingAR}
\end{figure}
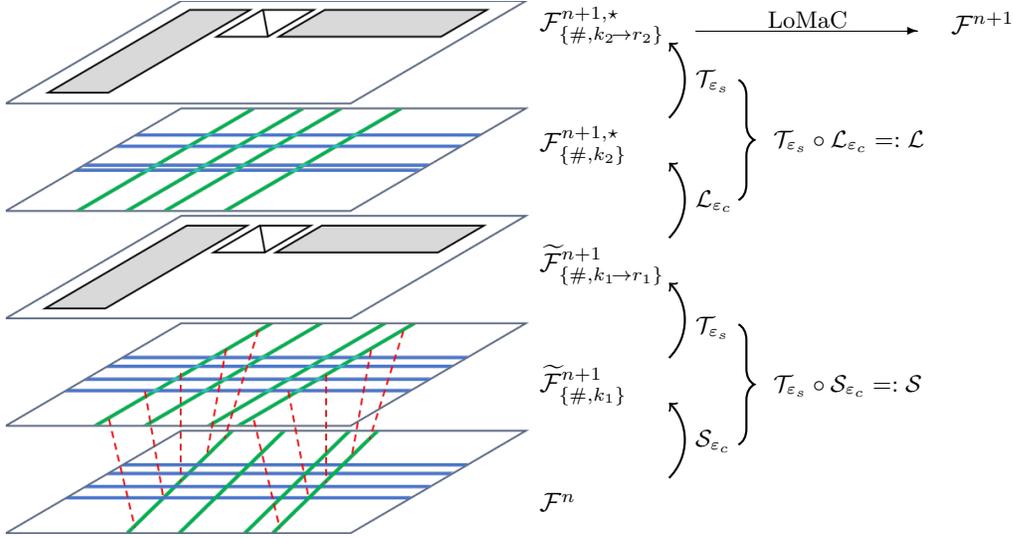

\Cref{fig:SamplingAR} illustrates the full adaptive-rank update procedure at each BE step. From bottom to top, the five layers are interpreted as follows. The process begins with the solution \( \mathcal{F}^n \) at the current time level. The operator \( \mathcal{S}_{\varepsilon_c} \) applies the ACA algorithm with sampling threshold \( \varepsilon_c \), performing local SL characteristic tracing and spatial reconstruction at the feet of the characteristics. This yields a CUR approximation of rank \( k_1 \), denoted \( \widetilde{\mathcal{F}}^{n+1}_{\{\#,k_1\}} \). Applying SVD truncation with tolerance \( \varepsilon_s \) via \( \mathcal{T}_{\varepsilon_s} \), we obtain a rank-\( r_1 \) approximation \( \widetilde{\mathcal{F}}^{n+1}_{\{\#,k_1 \to r_1\}} \). Next, the operator \( \mathcal{L}_{\varepsilon_c} \) performs ACA again, now acting entry-wise on the nonlinear update operator, as defined in \eqref{eq:BE_final}, producing a rank-\( k_2 \) approximation \( \mathcal{F}^{n+1,\star}_{\{\#,k_2\}} \). Another SVD truncation reduces it to rank \( r_2 \), yielding \( \mathcal{F}^{n+1,\star}_{\{\#,k_2 \to r_2\}} \). Finally, the LoMaC correction step is applied to restore conservation properties, resulting in the final low-rank solution \( \mathcal{F}^{n+1} \). 

The ACA+SVD strategy serves as a general compression framework that extends beyond the current setting. It is applicable to a broad class of problems where operators, solution updates, or nonlinear terms can be accessed locally without requiring full matrix assembly. This includes constructing adaptive-rank solvers for local time-stepping schemes, and developing compressed representations for nonlinear update formulas such as \eqref{eq:BE_final}. In particular, when ACA+SVD is applied to a local SL solver, the resulting operator reduces to the SLAR formulation proposed in our earlier work. The operator \( \mathcal{S} \) in this paper can thus be regarded as an instance of the SLAR method.

\begin{rem}
The ACA and SVD truncation tolerances directly determine the approximation rank and thus balance accuracy against computational and storage costs. A natural trade-off arises between discretization error and low-rank truncation error. For intrinsically low-rank solutions (or those well approximated in low-rank form), the error is dominated by discretization, so excessively small truncation tolerances are unnecessary. For solutions that are not inherently low-rank, the truncation tolerance should be chosen in line with the discretization error (or an estimate thereof), while simultaneously enforcing a maximum rank to limit memory consumption.
\end{rem}

\subsection{Preservation of Macroscopic Conservation Laws}
\label{subsec:moment preservation}

In our previous work \cite{zheng2025semi}, we ensured mass conservation in the SLAR method by simultaneously solving a macroscopic equation for the charge density, in a manner inspired by the LoMaC framework \cite{guo2024local}. A notable advantage of the SL method is its ability to take large time steps, which we preserved by discretizing the macroscopic equation with an implicit method. This resulted in a linear system for the charge density whose solution was used to correct the total mass of the distribution. In this work, we further extend this technique to simultaneously achieve mass, momentum, and energy conservation by implicitly solving a nonlinear system. In the following, we propose a nonlinear implicit solver for the macroscopic system; the updated macroscopic observables serve as constraints for correcting the low-rank solution, thereby ensuring local conservation of mass, momentum, and energy.

To demonstrate the idea, we consider a first-order BE time discretization, which offers a natural extension to high-order using DIRK methods. The process starts from a provisional non-conservative kinetic solution $f^{n+1,\star}$ obtained with the BE method \eqref{eq:BE_final}, which will be corrected to obtain a conservative solution. Motivated by property \eqref{eq:BGK invariants}, we seek a \textit{conservative kinetic solution} with the following LoMaC correction format:
\begin{equation}
    \label{eq:conservative kinetic solution format}
    f^{n+1}\left(\mathbf{U}^{n+1}\right) = f^{n+1,\star} - \mathcal{M}_{\mathbf{U}(f^{n+1,\star})} + \mathcal{M}_{\mathbf{U}^{n+1}}.
\end{equation}
Here, $\mathcal{M}_{\mathbf{U}^{n+1}}$ is the local Maxwellian constructed from $\mathbf{U}^{n+1}$, where $\mathbf{U}^{n+1}$ is obtained by solving the following macroscopic fluid system with a BE discretization
\begin{equation}
    \label{eq:BE macroscopic system}
    \mathbf{U}^{n+1} =  \mathbf{U}^n - \Delta t \mathbf{F}(\mathbf{U}^{n+1})_x,
\end{equation}
whose flux is \textit{implicitly} defined using the moments of the kinetic solution \eqref{eq:conservative kinetic solution format}, namely
\begin{equation}
    \label{eq:BE macroscopic flux}
    \mathbf{F}(\mathbf{U}^{n+1}) :=\left\langle\left(\begin{array}{l}
        f^{n+1} v\\
        f^{n+1} v^2\\
        f^{n+1}\left(\frac{1}{2} v^3\right)
    \end{array}\right)\right\rangle_v.
\end{equation}
Using the definition \eqref{eq:conservative kinetic solution format}, we can equivalently express the flux \eqref{eq:BE macroscopic flux} as
\begin{equation}
\label{eq:BE macroscopic flux expanded}
\begin{split}
    \mathbf{F}(\mathbf{U}^{n+1}) &:= \left\langle \left[f^{n+1,\star}-\mathcal{M}_{\mathbf{U}(f^{n+1,\star})}\right]
    \begin{pmatrix}
        0 \\
        0 \\
        \frac{1}{2} v^{3}
    \end{pmatrix}
    \right\rangle_{v}+\left\langle\mathcal{M}_{\mathbf{U}^{n+1}}\begin{pmatrix}
         v\\
        v^2\\
        \frac{1}{2} v^3
    \end{pmatrix}\right\rangle_v.
\end{split} 
\end{equation}
It is important to note that the first term on the right-hand side of \eqref{eq:BE macroscopic flux expanded} is fixed with respect to $\mathbf{U}^{n+1}$. In addition, the nonlinear dependence on $\mathbf{U}^{n+1}$ is entirely captured in the second term, which means that the non-conservative kinetic solution $f^{n+1,\star}$ is only used to initialize the closure of the fluid system. An iterative method then dynamically updates the kinetic solution through the LoMaC correction formula \eqref{eq:conservative kinetic solution format}, along with the macroscopic flux \eqref{eq:BE macroscopic flux expanded}, until nonlinear convergence is achieved on the system \eqref{eq:BE macroscopic system}. \Cref{fig:Fluid correction diagram} illustrates one step of the proposed method with the BE discretization, including the LoMaC correction loop, which is shown using solid lines. Note that the time level has been suppressed for brevity. Here dashed lines indicate that the kinetic solver is isolated from this process. Further details regarding the numerical treatment of this nonlinear system are provided in the following subsections.

\begin{figure}[!h]
\begin{center}
\begin{tikzpicture}

    \draw (2,0) rectangle (6,1); 
    \node at (2+2,0.5) {$f_t+vf_x=\frac{\mathcal{M}_{\mathbf{U}(f)}-f}{\epsilon}$};

    \draw (3.4,2) rectangle (4.6,3); 
    \node at (2+2,2.5) {$f^{\star}$}; 

    \draw (8.75,0) rectangle (13.0,1); 
    \node at (10.9,0.5) {$\mathbf{U}_{t}+\mathbf{F}(\mathbf{U}(f))_x=0$}; 

    \draw (8.0,2) rectangle (14.25,3); 
    \node at (11,2.5) {$f(\mathbf{U}) = f^{\star}-\mathcal{M}_{\mathbf{U}(f^{\star})}+\mathcal{M}_{\mathbf{U}}$}; 

    \node at (11,-1.5) {\textbf{Converged?}}; 
    \node at (14.8,-1.8) {\scriptsize No};
    \node at (7.0,-1.8) {\scriptsize Yes};
    
    \node[anchor=west] at (4,1.5) {\scriptsize Adaptive-rank non-conservative update \eqref{eq:BE_final}};
    \node[anchor=west] at (11,-0.65) {\scriptsize Conservative update for the system \eqref{eq:BE macroscopic system}};
    \node[anchor=west] at (11,1.5) {\scriptsize Update macroscopic flux \eqref{eq:BE macroscopic flux expanded}};

    \node[anchor=west] at (14.3,2.8) {\scriptsize Modify closure \eqref{eq:conservative kinetic solution format}};

    \node[anchor=west] at (4,-0.75) {\scriptsize Prepare for the next time step
};
    \node at (6.2, 2.8) {\scriptsize Initialize closure \eqref{eq:conservative kinetic solution format}};

    \draw[->][thick, dashed] (4,1) -- (4,2);

    \draw[->][thick, dashed] (4.6,2.5) -- (8.0,2.5);

    \draw[->][thick] (11,2) -- (11,1);

    \draw[->][thick] (11,0) -- (11,-1.25);

    \draw[thick] (12.5,-1.5) -- (16+1,-1.5);
    \draw[thick] (17,-1.5) -- (16+1,2.5);
    \draw[->][thick] (17,2.5) -- (14.25,2.5);

    \draw[thick, dashed] (9.5,-1.5) -- (4,-1.5);
    \draw[->][thick, dashed] (4,-1.5) -- (4,0);

\end{tikzpicture}
\end{center}
\caption{Algorithm flowchart for the BGK equation with a BE discretization and LoMaC correction.}
\label{fig:Fluid correction diagram}
\end{figure}

\subsubsection{Spatial Discretization of the Nonlinear system}

There are many conservative spatial discretizations available for approximating the macroscopic system. To obtain a fully discrete method for \eqref{eq:BE macroscopic system} based on the BE method, we combine an upwind flux-splitting method with a high-order reconstruction technique. This combination produces a high-order conservative approximation to the spatial derivative, which we write as
\begin{equation}
    \label{eq:BE macroscopic fully discrete}
    \mathbf{U}_i^{n+1} = \mathbf{U}_i^n - \frac{\Delta t}{\Delta x}\left(\widehat{\mathbf{F}}^{n+1}_{i+\frac{1}{2}}-\widehat{\mathbf{F}}^{n+1}_{i-\frac{1}{2}}\right).
\end{equation}
Here, the numerical fluxes $\widehat{\mathbf{F}}^{n+1}_{i\pm\frac{1}{2}}$ for the macroscopic system are defined using a kinetic flux vector splitting (KFVS) technique \cite{1994Kinetic,1995Gas}, which separates the velocity domain into positive and negative components, corresponding to particle motion in different directions. This provides the correct wind directions needed to perform high-order spatial reconstruction methods. In KFVS, we express the flux as $\widehat{\mathbf{F}}^{n+1} = \widehat{\mathbf{F}}^{n+1,+} + \widehat{\mathbf{F}}^{n+1,-}$, where the split fluxes used in the reconstructions are given by
\begin{equation}
    \label{eq:macroscopic flux stencil points}
    \widehat{\mathbf{F}}^{n+1,\pm} = \widetilde{\int_{\mathbb{R}^{\pm}}} f^{n+1}\left(\begin{array}{c}
             v\vspace{-0.2cm}  \\
             v^2\vspace{-0.2cm}\\
             \frac{1}{2}v^3
        \end{array}\right)\,dv.
\end{equation}
Here, $\widetilde{\int}_{\mathbb{R}^{\pm}} \cdot\,dv$ denotes a low-rank numerical integration over the half-real velocity domains, as described in \Cref{app:low_rank_integration}. Once the split flux stencil values are obtained from \eqref{eq:macroscopic flux stencil points}, a high-order reconstruction is applied to evaluate the fluxes at cell interfaces. In this work, we use the fifth-order classical WENO-JS reconstructions \cite{JIANG1996202}, which are written as
\begin{equation}
    \begin{split}
        \label{eq:macroscopic flux splitting}
        \widehat{\mathbf{F}}^{n+1}_{i+\frac{1}{2}} &= \widehat{\mathbf{F}}^{n+1,+}_{i+\frac{1}{2}} + \widehat{\mathbf{F}}^{n+1,-}_{i+\frac{1}{2}}\\
        &:=\mathcal{R}^{+}\left(\left\{\mathbf{F}^{n+1,+}_{i-2}, \cdots, \mathbf{F}^{n+1,+}_{i+2}\right\} \right) + \mathcal{R}^{-}\left(\left\{\mathbf{F}^{n+1,-}_{i-1}, \cdots, \mathbf{F}^{n+1,-}_{i+3}\right\}\right).
    \end{split}
\end{equation}
Here, we use $\mathcal{R}^{\pm}$ to denote the WENO reconstructions of the flux at a cell interface. 

Spatial discretization of \eqref{eq:BE macroscopic fully discrete} results in a system of nonlinear algebraic equations for the discrete macroscopic variables $\mathbf{U}^{n+1}$, which can be written in the form
\begin{equation}
    \label{eq:nonlinear system BE macroscopic}
    \mathbf{G}\Big( \mathbf{U}^{n+1} \Big) = 0.
\end{equation}
An important observation is that the size of the resulting nonlinear problem scales with the size of the system for the macroscopic variables and \textit{not} the kinetic equation. In other words, it is of reduced dimensionality.

\subsubsection{Newton-Krylov Methods for the Nonlinear System}
\label{subsec:JFNK}

In the previous section, we found that the implicit, fully-discrete macroscopic system for the BGK equation with the dynamic closure leads to a nonlinear system of the form
\begin{equation}
    \label{eq:General nonlinear problem}
    \mathbf{G}\Big( \mathbf{U} \Big) = 0,
\end{equation}
where we have dropped the reference to the time level, for brevity. A Newton method for \eqref{eq:General nonlinear problem} uses the linearization about a current iterate $\mathbf{U}^{(k)}$, resulting in the sequence of linear systems
\begin{equation}
    \label{eq:Newton linear system}
    J^{(k)} \delta \mathbf{U}^{(k)} = -\mathbf{G}\Big( \mathbf{U}^{(k)} \Big), \quad \text{for } k = 1, 2, \cdots,
\end{equation}
where $J^{(k)} = \mathbf{G}'\Big( \mathbf{U}^{(k)} \Big)$ is the Jacobian matrix of the system evaluated at the current iterate and $\delta \mathbf{U}^{(k)}$ is the correction to $\mathbf{U}^{(k)}$. After solving \eqref{eq:Newton linear system} for the correction $\delta \mathbf{U}^{(k)}$, we obtain the updated state
\begin{equation}
    \label{eq:Newton solution update}
    \mathbf{U}^{(k+1)} = \mathbf{U}^{(k)} + \delta \mathbf{U}^{(k)}.
\end{equation}

As discussed at the end of the previous section, the size of this nonlinear system depends entirely on the size of the problem for the macroscopic system. In the more general setting, the system will be of size $M_{d}N^{d} \times M_{d}N^{d}$, where $N$ is the number of cells per spatial dimension, $d=1,2,3$, and $M_d$ is the number of conserved quantities with dimension $d$. Other complications arise due to the use of nonlinear spatial reconstruction methods in the macroscopic solver and its nonlinear coupling to the solution from the kinetic solver. These features make analytical evaluations of the Jacobian expensive.

Instead, we consider Newton-Krylov methods, which typically combine an inexact Newton method with Krylov subspace methods. Inexact Newton methods use an approximation of the Jacobian, which slightly weakens its convergence rate compared to an exact Newton method. However, the adaptive-rank kinetic solver we employ provides a high-quality, inexpensive initial guess for the Newton method that is non-conservative. A notable advantage of inexact Newton methods is that they can be used in an entirely matrix-free setting and are known as Jacobian-free Newton-Krylov methods \cite{knoll2004jacobian,knoll2005jacobian}. In particular, products between the Jacobian matrix $J^{(k)}$ and an arbitrary vector $\mathbf{V}$ can be approximated as
\begin{equation}
    \label{eq:Jacobian-vector product}
    J^{(k)} \mathbf{V} = \frac{\mathbf{G}\Big( \mathbf{U}^{(k)} + \epsilon \mathbf{V} \Big) - \mathbf{G}\Big( \mathbf{U}^{(k)} \Big)}{\epsilon} + \mathcal{O}(\epsilon),
\end{equation}
without ever explicitly forming a Jacobian matrix. Here the parameter $\epsilon > 0$ (different from the Knudsen number) defines a perturbation of the current state $\mathbf{U}^{(k)}$, with the error being proportional to $\epsilon$. This approximation avoids the explicit formation of the Jacobian matrix. Consequently, it naturally motivates the use of Krylov subspace methods, in particular GMRES, which search for a solution to the linear systems \eqref{eq:Newton linear system} using the subspaces spanned by recursive Jacobian-vector products.

\begin{rem}
In our implementation, the nonlinear tolerance of the Newton–Krylov solver is set to machine precision to ensure exact conservation of the moments. The performance is largely insensitive to the tolerance of the inner Krylov solver. While previous studies have considered dynamically tightening this tolerance during the nonlinear iterations to avoid oversolving a Newton step \cite{knoll2004jacobian}, our experiments showed no practical benefit from such strategies. A distinguishing feature of our approach is that the Newton iteration is initialized with a high-quality guess provided by the low-rank solver, which we found to be the primary factor in reducing the number of iterations.
\end{rem}

\subsection{Extensions to High-order using Runge-Kutta Methods}
\label{subsec:high-order DIRK}

To obtain higher-order temporal accuracy, we apply high-order DIRK methods \cite{alexander1977diagonally,kennedy2016diagonally}. The Butcher tableau for a DIRK method is given by
$
\begin{array}{c|c}
\mathbf{c} & A \\
\hline
 & \mathbf{b}^T
\end{array},
$
where $A=(a_{i,j})\in\mathbb{R}^{s\times s}$ is an invertible lower triangular matrix, $\mathbf{c}=[c_1,\ldots,c_s]^\top$ is the intermediate coefficient vector, and $\mathbf{b}^\top=[b_1,\ldots,b_s]$ represents the quadrature weights. In order to address the stiffness in \eqref{eq:BGK} when $\epsilon \ll 1$, we use DIRK methods with the SA property, for which $A(s,:) = \mathbf{b}^\top$ and $f^{n+1} = f^{(s)}$.

The SL-FD method from our previous work \cite{li2023high} applies DIRK methods along the characteristics, so that the solution at internal stages $f^{(k)}$ is given by
\begin{align}
    \label{eq:SLFD_DIRK_form1}
    \begin{split}
    f^{(k)}_{i,j} &= f^n(x_i-v_jc_k\Delta t,v_j) + \Delta t\sum\limits_{\ell=1}^{k}a_{k,\ell}\left(\frac{\mathcal{M}_{\mathbf{U}({f^{(\ell)}})}-f^{(\ell)}}{\epsilon}\right)(x_i-v_j(c_k-c_\ell)\Delta t,v_j) \\
    &=:\widetilde{f}^{(k)} + a_{k,k}\Delta t\left(\frac{\mathcal{M}_{\mathbf{U}({f^{(k)}})}-f^{(k)}}{\epsilon}\right)(x_i,v_j),\quad \text{for}~~ k = 1,\ldots,s.
    \end{split}
\end{align}
The implicit advance \eqref{eq:SLFD_DIRK_form1} at each stage can be written in an explicit form, similar to the BE discretization \eqref{eq:BE_final}, which leads to the stage-wise update
\begin{equation}\label{eq:SLFD_DIRK_final}
    f^{(k)}_{i,j} = \frac{\epsilon \widetilde{f}^{(k)}_{i,j} + a_{k,k}\Delta t \left(\mathcal{M}_{\mathbf{U}\left(\widetilde{f}^{(k)}\right)}\right)_{i,j}}{\epsilon + a_{k,k} \Delta t}, \quad \text{for } k = 1, \ldots, s.
\end{equation}
As a consequence of the SA property, it follows that $f^{n+1} = f^{(s)}$. Additionally, we note that there is an alternative approach to evaluating the intermediate solutions $\widetilde{f}^{(k)}$ in \eqref{eq:SLFD_DIRK_final}, known as the Shu-Osher form \cite{ding2023accuracy}. In this study, we adopt the Shu-Osher form as it proved more reliable for simulations involving discontinuous solutions.

High-order time discretizations for the macroscopic system can be achieved in a similar manner by extending the definitions \eqref{eq:conservative kinetic solution format}–\eqref{eq:BE macroscopic flux} to internal stages of the DIRK methods. At each stage, the intermediate solution $f^{(k),\star}$ is computed using the adaptive-rank method with greedy sampling, as described in \Cref{subsec:SLAR for the BGK equation}, where the same sampling and truncation strategy is applied independently at each stage.

The full adaptive-rank update procedure within each DIRK stage consists of four successive low-rank truncations followed by a final macroscopic correction. The stage-$k$ procedure is summarized as:
\begin{equation*}
    \left\{\mathcal{F}^{n}, \mathcal{F}^{(1)}, \cdots, \mathcal{F}^{(k-1)} \right\}
    \xrightarrow{\;\mathcal{S}_{\epsilon_{c}}\;} 
    \widetilde{\mathcal{F}}^{(k)}_{\{\#,k_{1}\}}
    \xrightarrow{\;\mathcal{T}_{\epsilon_{s}}\;} 
    \widetilde{\mathcal{F}}^{(k)}_{\{\#,k_{1} \to r_{1}\}}
    \xrightarrow{\;\mathcal{L}_{\epsilon_{c}}\;} 
    \mathcal{F}^{(k),\star}_{\{\#,k_{2}\}}
    \xrightarrow{\;\mathcal{T}_{\epsilon_{s}}\;} 
    \mathcal{F}^{(k),\star}_{\{\#,k_2 \to r_2\}}
    \xrightarrow{\;\text{LoMaC}\;} 
    \mathcal{F}^{(k)}.
\end{equation*}
Here, \( \left\{\mathcal{F}^{n}, \mathcal{F}^{(1)}, \cdots, \mathcal{F}^{(k-1)} \right\} \) denotes the set of available low-rank numerical solutions from previous stages. Only the final corrected solution $\mathcal{F}^{(k)}$ is passed to the next stage.

For the LoMaC correction, the conservative kinetic solution at each stage is defined by
\begin{equation}
    \label{eq:conservative kinetic solution DIRK stage format}
    \mathcal{F}^{(k)} = \mathcal{F}^{(k),\star} - \mathcal{M}_{\mathbf{U}(\mathcal{F}^{(k),\star})} + \mathcal{M}_{\mathbf{U}^{(k)}}, \quad\text{for }k=1,\ldots,s,
\end{equation}
where the macroscopic solution $\mathbf{U}^{(k)}$ at the internal stages is calculated according to
\begin{equation}
    \label{eq:Macroscopic DIRK}
    \mathbf{U}^{(k)} = \mathbf{U}^{n} - \Delta t\sum_{\ell = 1}^{k} a_{k,\ell} \mathbf{F}(\mathbf{U}^{(\ell)})_x, \quad\text{for }k=1,\ldots,s,
\end{equation}
and the macroscopic fluxes at the stages are defined as
\begin{equation}
    \label{eq:DIRK macroscopic flux}
    \mathbf{F}(\mathbf{U}^{(k)}) := \left\langle \mathcal{F}^{(k)}
    \begin{pmatrix}
        v \\
        v^2 \\
        \frac{1}{2}v^{3}
    \end{pmatrix}
    \right\rangle_{v}, \quad\text{for }k=1,\ldots,s.
\end{equation}

We also note that there are some structural differences in the nonlinear system that arise when using DIRK methods. In the case of the fully discrete BE discretization, we need to solve a nonlinear problem \eqref{eq:nonlinear system BE macroscopic} at the new time level. Correspondingly, for the high-order DIRK methods, we need to solve a sequence of nonlinear problems at each stage of the method, which we denote as
\begin{equation*}
    \mathbf{G}^{(k)}\Big( \mathbf{U}^{(k)} \Big) = 0, \quad\text{for }k=1,\ldots,s.
\end{equation*}
Here $\mathbf{G}^{(k)}$ is used to indicate that a distinct nonlinear system must be solved at each DIRK stage. These problems can be solved using the same methods discussed in \Cref{subsec:JFNK}.

\subsection{Asymptotic Analysis}
\label{subsec:asymptotic analysis}

In this section, we analyze the properties of the proposed adaptive-rank method in the fluid limit. First, we recall key asymptotic properties of the full-grid high-order local SL method, as established in our previous work \cite{ding2023accuracy}, which serves as the starting point for this study. We then show that the adaptive closure technique, used to preserve the conservation laws, also preserves the asymptotic properties of the high-order local SL method in the full-grid setting. We also analyze the impact of low-rank approximation on the asymptotic behavior. By utilizing error estimates for low-rank matrix approximations, we derive corresponding estimates that establish a relationship between the Knudsen number $\epsilon$ and the truncation errors of the CUR and SVD approximations. The local Maxwellian which characterizes the fluid limit can not be represented in terms of separable functions, except in limited circumstances. From these results, a conditional AP property of the low-rank scheme is established, and we suggest an alternative which removes these restrictions.

We first recall the AA property of the temporal discretization used in the local SL solver. As in \cite{ding2023accuracy} (see Theorems 3.5 and 3.6 therein), we analyze the scheme under the assumption of exact spatial discretization. Specifically, we assume that the phase space mesh has a tensor-product structure and that there is no spatial interpolation error in the reconstruction at the feet of the characteristics. The AA property holds under the additional assumption that the initial data is well-prepared, which we define as follows.
\begin{defn}
Let $\epsilon > 0$ denote a small parameter. The initial data $f_{0}^{\epsilon}$ are said to be well-prepared if they admit the formal asymptotic expansion
\begin{equation*}
    f_{0}^{\epsilon} = f_{0} + \epsilon f_{1} + \epsilon^2 f_{2} + \cdots
\end{equation*}
where each $f_{j}$ is consistent with the order of the corresponding term in the expansion. The assumption of well-prepared initial data prevents the formation of spurious initial layers or fast transients in the limit $\epsilon \rightarrow 0$ that would otherwise need to be resolved by the discretization.
\end{defn}

Next, we formally state the AA property of the temporal discretization for the full-grid kinetic scheme, which is as follows.
\begin{thm}
    \label{thm:AA property for the local solver}
    Suppose that the initial data are well-prepared. Consider the kinetic SL scheme defined in \eqref{eq:SLFD_DIRK_form1}, using an $s$-stage SA DIRK method, and assume that exact spatial interpolation is used at the feet of the characteristics. Then, in the limit $\epsilon \to 0$, the temporal accuracy of the kinetic scheme implies the corresponding temporal accuracy of the macroscopic scheme defined by \eqref{eq:Macroscopic DIRK} and \eqref{eq:DIRK macroscopic flux}, provided that the Butcher tableau satisfies the following order conditions (valid up to third order):
    \begin{align*}
        &\text{First-order: } c_{s} = 1, \\
        &\text{Second-order: } d_{s} = \frac{1}{2}, \\
        &\text{Third-order: } g_{s} = h_{s} = \frac{1}{6}, \quad G_{s} = \frac{1}{6},
    \end{align*}
    where the coefficients $c_{s}$, $d_{s}$, $g_{s}$, $h_{s}$, and $G_{s}$ are defined recursively according to equations (3.10) and (3.24) in \cite{ding2023accuracy}.
\end{thm}

A unique component of this work is its use of the macroscopic system to provide a self-consistent correction to the moments of the kinetic solution. This process enables the simultaneous preservation of mass, momentum, and energy up to a user-specified nonlinear tolerance. As discussed in \Cref{subsec:moment preservation}, the high-order local SL method provides a high-quality initial guess, and a ``small" correction is applied to preserve the moments. The next result shows that this process preserves the asymptotic properties of the high-order local SL method in the full-grid setting.
\begin{thm}
\label{thm:moment preservation AA property}

Suppose the initial data are well-prepared, and assume that the fully-discrete kinetic scheme is AP but does not simultaneously preserve mass, momentum, and energy. Then, the self-consistent macroscopic correction procedure defined in \eqref{eq:Macroscopic DIRK} and \eqref{eq:DIRK macroscopic flux}, when applied to the intermediate solutions $f^{(k)}$ from \eqref{eq:conservative kinetic solution DIRK stage format}, preserves the underlying AP property of the scheme. Moreover, if the full-discrete kinetic scheme is also AA, then the macroscopic correction preserves the AA property as well.


\begin{proof}
    
    Suppose that the fully-discrete scheme is AP. Then, by definition, it follows that the discrete kinetic solution satisfies
    \begin{equation}
        \label{eq:AA property of the local solver}
        F^{(k),\star} = \mathcal{M}_{\mathbf{U}\left( F^{(k),\star} \right)} + \mathcal{O}(\epsilon), \quad\text{for }k=1,\ldots,s.
    \end{equation}
    Here, we use $F^{(k),\star}$ to distinguish the non-conservative kinetic solution from its fully-conservative (mass, momentum, and energy) counterpart. 

    Now consider the self-consistent macroscopic stage values $\mathbf{U}^{(k)}$ defined by the moment equations \eqref{eq:Macroscopic DIRK} and \eqref{eq:DIRK macroscopic flux}. These are closed using the corrected kinetic distribution $F^{(k)}$, which is constructed following the conservative format \eqref{eq:conservative kinetic solution DIRK stage format} and is given by
    \begin{equation*}
        F^{(k)}(\mathbf{U}^{(k)}) = F^{(k),\star} - \mathcal{M}_{\mathbf{U}\left( F^{(k),\star} \right)} + \mathcal{M}_{\mathbf{U}\left( F^{(k)} \right)}.
    \end{equation*}
    Substituting \eqref{eq:AA property of the local solver} into this expression and taking the limit $\epsilon \to 0$, we obtain
    \begin{equation*}
            F^{(k)}\left(\mathbf{U}^{(k)}\right) = \mathcal{M}_{\mathbf{U}^{(k)}}, \quad\text{for }k=1,\ldots,s.
    \end{equation*}
    This limiting solution provides the necessary closure to the macroscopic flux defined by \eqref{eq:DIRK macroscopic flux}, namely
    \begin{equation*}
        \mathbf{F}(\mathbf{U}^{(k)}) := \left\langle \mathcal{M}_{\mathbf{U}^{(k)}}
        \begin{pmatrix}
            v \\
            v^2 \\
            \frac{1}{2}v^{3}
        \end{pmatrix}
        \right\rangle_{v}, \quad\text{for }k=1,\ldots,s.
    \end{equation*}
    When combined with the self-consistent discretization of the system \eqref{eq:Macroscopic DIRK}, we obtain the compressible Euler equations. Therefore, the correction process preserves the AP property of the kinetic scheme. If the fully-discrete kinetic scheme is also AA, then the discretization of the self-consistent macroscopic system converges at the same order of accuracy as the kinetic scheme, since the same DIRK method is employed for both systems. Furthermore, the deviations from the local Maxwellian are again $\mathcal{O}(\epsilon)$, so this property remains unaffected.
    
\end{proof}

\end{thm}


In the low-rank setting, the use of low-rank approximations may cause the solution to deviate from the correct asymptotic limit. Each of the approximations contributes to the overall truncation error, which we now quantify in the following result.

\begin{thm}
\label{thm:Error between F and Maxwellian}
Suppose that a SA DIRK temporal discretization with the AA property, satisfying the conditions of \Cref{thm:AA property for the local solver}, is paired together with a local SL solver. Then, at the end of each stage $k$, the conservative low-rank solution \(\mathcal{F}^{(k)}\) defined following \eqref{eq:conservative kinetic solution DIRK stage format}, satisfies the following truncation error bound
\begin{equation*}
\norm{ \mathcal{F}^{(k)} - \mathcal{M}_{\mathbf{U}(\mathcal{F}^{(k)})} } \lesssim \epsilon + \norm{ E_{\epsilon_{s}} \left( \mathcal{F}_{\{ \#,k_{2} \}}^{(k)} \right) } + \norm{ E_{\epsilon_{c}} \left( F^{(k)} \right) },
\end{equation*}
where we denote the truncation errors due to the CUR and SVD approximations as
\begin{equation*}
    E_{\epsilon_{s}} \left( \mathcal{F}_{\{ \#,k_{2} \}}^{(k)} \right) := \mathcal{F}_{\{ \#,k_{2} \to r_{2} \}}^{(k)} - \mathcal{F}_{\{ \#,k_{2} \}}^{(k)}, \quad E_{\epsilon_{c}} \left( F^{(k)} \right) := \mathcal{F}_{\{ \#,k_{2} \}}^{(k)} - F^{(k)}.
\end{equation*}

\begin{proof}
Using the definition of the conservative form of the low-rank solution \eqref{eq:conservative kinetic solution DIRK stage format}, we can decompose the error as
\begin{align}
    \mathcal{F}^{(k)} 
  - \mathcal{M}_{\mathbf{U}(\mathcal{F}^{(k)})} &= \mathcal{F}_{\{ \#,k_{2} \to r_{2} \}}^{(k),\star} 
  - \mathcal{M}_{\mathbf{U}\left( \mathcal{F}_{\{ \#,k_{2} \to r_{2} \}}^{(k),\star} \right)}, \nonumber \\
  &= F^{(k)} - \mathcal{M}_{\mathbf{U}\left( F^{(k)} \right)} + \mathcal{F}_{\{ \#,k_{2} \to r_{2} \}}^{(k),\star} 
  - F^{(k)} + \mathcal{M}_{\mathbf{U}\left( F^{(k)} \right)} -  \mathcal{M}_{\mathbf{U}\left( \mathcal{F}_{\{ \#,k_{2} \to r_{2} \}}^{(k),\star} \right)}. \label{eq:stage-wise difference between F and M_U}
\end{align}

Since the initial data is well-prepared and the local solver is AA, then we immediately have that
\begin{equation}
    \label{eq:AA property of the local solver discrete}
    \norm{ F^{(k)} - \mathcal{M}_{\mathbf{U}\left( F^{(k)} \right)} } \lesssim \epsilon.
\end{equation}

The second term can be estimated by further decomposing the error into the two parts. One concerns the SVD truncation of a CUR object, while the other represents the approximation error for the CUR decomposition of full-rank data. This gives
\begin{align}
    \mathcal{F}_{\{ \#,k_{2} \to r_{2} \}}^{(k),\star} - F^{(k)} &= \mathcal{F}_{\{ \#,k_{2} \to r_{2} \}}^{(k),\star} - \mathcal{F}_{\{ \#,k_{2} \}}^{(k),\star} + \mathcal{F}_{\{ \#,k_{2} \}}^{(k),\star}   - F^{(k)} := E_{\epsilon_{s}} \left( \mathcal{F}_{\{ \#,k_{2} \}}^{(k),\star} \right) + E_{\epsilon_{c}} \left(  F^{(k)} \right). \label{eq:final layer error SVD of CUR versus true data}
\end{align}

The last term involves a difference between two Maxwellians, which can be estimated using the Mean Value Theorem and the error equation \eqref{eq:final layer error SVD of CUR versus true data}, together with the triangle inequality. That is, there exists an $\mathcal{F}_{*}^{(k)}$ such that
\begin{align}
\label{eq:Bound of differences between Maxwellians}
    \norm{\mathcal{M}_{\mathbf{U}\left(F^{(k)}\right)} - \mathcal{M}_{\mathbf{U}\left( \mathcal{F}^{(k)}_{\{\#,k_{2} \to r_{2}\}} \right)}}  &\leq \norm{\mathcal{M}_{\mathbf{U}\left( \mathcal{F}_{*}^{(k)} \right)}^{'} } \left( \norm{ E_{\epsilon_{s}} \left( \mathcal{F}_{\{ \#,k_{2} \}}^{(k),\star} \right) } + \norm{ E_{\epsilon_{c}} \left(  F^{(k)} \right) } \right).
\end{align}
In the above, $\mathcal{M}_{\mathbf{U}\left( \mathcal{F}_{*}^{(k)} \right)}^{'}$ denotes the Jacobian matrix of a local Maxwellian, which will be bounded provided that the moments of $\mathcal{F}_{*}^{(k)}$ satisfy $\rho_{*}^{(k)} \geq \rho_{m} > 0$ and $T_{*}^{(k)} \geq T_{m} > 0$.

Finally, applying the triangle inequality to \eqref{eq:stage-wise difference between F and M_U} together with each of the inequalities in \eqref{eq:AA property of the local solver discrete} to \eqref{eq:Bound of differences between Maxwellians}, we obtain the result provided in the statement. This completes the proof.
\end{proof}

\end{thm}

The term $E_{\epsilon_{s}} \left( \mathcal{F}_{\{ \#,k_{2} \}}^{(k),\star} \right)$ in \eqref{eq:final layer error SVD of CUR versus true data} requires some additional explanation. Recall \eqref{eq:SLFD_DIRK_final}, which provides the following element-wise definition to build the CUR decomposition $\mathcal{F}_{\{ \#,k_{2}\}}^{(k),\star}$. In the limit $\epsilon \to 0$, this formula reduces to
\begin{equation*}
    \left( \mathcal{F}_{\{ \#,k_{2}\}}^{(k),\star} \right)_{i,j} = \left( \mathcal{M}_{\mathbf{U}\left( \widetilde{\mathcal{F}}_{\{ \#,k_{1} \to r_{1}\}}^{(k)} \right)} \right)_{i,j}.
\end{equation*}
As $\widetilde{\mathcal{F}}_{\{ \#,k_{1} \to r_{1}\}}^{(k)}$ is constructed using both CUR and SVD approximations, there will be an accumulation of the truncation errors due to earlier stages, resulting in a perturbed local Maxwellian that will be approximated by a CUR decomposition. Due to the interpolation property of the CUR decomposition, the previous limiting equality will be enforced only at selected indices and will be approximate elsewhere. To analyze the truncation error at the interpolation points, let
\begin{align*}
    \widetilde{\mathcal{F}}_{\{ \#,k_{1} \to r_{1}\}}^{(k)} &= \widetilde{\mathcal{F}}_{\{ \#,k_{1} \to r_{1}\}}^{(k)} - \widetilde{\mathcal{F}}_{\{ \#,k_{1}\}}^{(k)} + \widetilde{\mathcal{F}}_{\{ \#,k_{1}\}}^{(k)} - \widetilde{F}^{(k)} + \widetilde{F}^{(k)}
    := E_{\epsilon_{s}}\left( \widetilde{\mathcal{F}}_{\{ \#,k_{1}\}}^{(k)} \right) + E_{\epsilon_{c}}\left( \widetilde{F}^{(k)} \right) + \widetilde{F}^{(k)}.
\end{align*}
Again, by the Mean Value Theorem, there exists (another) $\mathcal{F}_{*}^{(k)}$ such that
\begin{equation}
\label{eq:Linearization of a local Maxwellian}
    \left( \mathcal{F}_{\{ \#,k_{2}\}}^{(k),\star} \right)_{i,j} = \left( \mathcal{M}_{\mathbf{U}\left( \widetilde{F}^{(k)} \right)} \right)_{i,j} + \Bigg( \mathcal{M}_{\mathbf{U}\left( \mathcal{F}_{*}^{(k)} \right)}^{'} \left[ E_{\epsilon_{s}}\left( \widetilde{\mathcal{F}}_{\{ \#,k_{1}\}}^{(k)} \right) + E_{\epsilon_{c}}\left( \widetilde{F}^{(k)} \right) \right] \Bigg)_{i,j}.
\end{equation}
Although this equality holds only at select interpolation points, it illustrates the accumulation of truncation error from the earlier stages of the method. As the subsequent SVD truncation destroys the interpolation property of CUR, it is difficult to say more about this accumulation of error.

We now state the main limiting property of the proposed method which appears as a corollary of \Cref{thm:Error between F and Maxwellian}. More specifically, we establish a \textit{conditional} AA property, which is analytically satisfied if the limiting solution exhibits rank degeneracy.
\begin{cor}
\label{cor:conditional AA property}
Let an $s$-stage SA DIRK method with the AA property (satisfying the conditions of \Cref{thm:AA property for the local solver}) be coupled with a local SL solver and the macroscopic correction procedure defined by \eqref{eq:conservative kinetic solution DIRK stage format} to \eqref{eq:DIRK macroscopic flux}. Then the proposed low-rank scheme satisfies the AA property in the limit $\epsilon \to 0$, provided that for each stage $k = 1,\ldots,s$, the following condition holds:
\begin{equation*}
    \norm{ E_{\epsilon_{s}} \left( \mathcal{F}_{\{ \#,k_{2} \}}^{(k)} \right) } + \norm{ E_{\epsilon_{c}} \left( F^{(k)} \right) } \lesssim \epsilon.
\end{equation*}
\end{cor}

\begin{rem}
To preserve the AA property of the full-grid scheme, the low-rank truncation errors must not perturb the distribution significantly away from the local Maxwellian. This condition is automatically satisfied for \emph{any} distribution if the truncation tolerances $\epsilon_{c}$ and $\epsilon_{s}$ are chosen to be $\mathcal{O}(\epsilon)$. However, doing so may undermine the efficiency of the low-rank method, as it may require large ranks to achieve such high accuracy. The key advantage of the proposed scheme is recovered when the limiting solution exhibits \emph{rank degeneracy} when $\epsilon \to 0$ — that is, when the local Maxwellian can be exactly written using separable functions. In this setting, it is possible to choose tolerances of size $\mathcal{O}(\epsilon)$ while maintaining a low computational cost.
\end{rem}

\begin{rem}
To alleviate the restriction on the truncation tolerances in \Cref{cor:conditional AA property}, one could alternatively consider a \textit{multiplicative} decomposition of the distribution function, i.e., $f = \mathcal{M}_{\mathbf{U}(f)}g$, and apply the low-rank decomposition only to the function $g$, which is a rank-1 function when $\epsilon \to 0$. This idea has been considered in other low-rank methods \cite{baumann2024stable,einkemmer2021efficient} for the case of isothermal Maxwellians with and without drift. The use of a multiplicative decomposition requires additional modifications to the scheme presented in this work, so we do not consider this approach here and leave it to future work.
\end{rem}

\section{Numerical Tests}\label{sec:numerical_tests}

In this section, we evaluate the proposed low-rank method using benchmark problems for the BGK model. We study the accuracy of the method across a range of $\epsilon$. Additionally, we demonstrate the capability of the method to capture shocks and discontinuous solution structures by applying them to a classical Riemann problem. To demonstrate the capabilities of the method in addressing multi-scale problems, we also consider a more challenging mixed regime problem. 

In all experiments, we use a four-stage third-order DIRK method, which is asymptotically accurate for the BGK model, and the interval $[-10,10]$ is used for the velocity domain. Note that the Butcher table for the scheme can be found in Table B10 of \cite{ding2023accuracy}). With regard to the JFNK method, we used a tolerance of $10^{-14}$ for the Newton method and a tolerance of $10^{-6}$ for the inner Krylov solver. The time step in the simulations is calculated from the CFL as $\Delta t = \text{CFL} \Delta x / \lvert v_{max}\rvert$. In certain examples, we measure the absolute error in a conserved quantity $q(x,t)$, e.g., mass, momentum, and energy, as $\left\lvert Q(t) - Q(0) \right\rvert$, where we have defined
\begin{align*}
    Q(t) = \int_{\Omega_{x}} q(x,t) \,dx \approx \Delta x \sum_{i=1}^{N_{x}} q_{i}(t).
\end{align*}

\begin{exa}(Consistent initial data \cite{li2023high}) Consider the BGK model with the initial Maxwellian distribution:
\begin{equation}\label{eq:consistent_IC}
f(x,v,0) = \frac{\rho_0}{\sqrt{2\pi T_0}}\exp\left(-\frac{(v-u_0)^2}{2T_0}\right),\quad x\in[-1,1]
\end{equation}
with initial density $\rho_0(x) = 1$, initial temperature $T_0(x) = 1$, and initial mean velocity:
\begin{equation*}
    u_0(x) = \frac{1}{10}\left[\exp\left(-(10x-1)^2\right)-2\exp\left(-(10x+3)^2\right)\right].
\end{equation*}
We apply periodic boundary conditions and set the final time to $t = 0.04$. The relative truncation tolerances for the CUR and SVD approximations are set to $\epsilon_C = 10^{-9}$ and $\epsilon_S = 10^{-8}$, respectively.

In \Cref{fig:consistent temporal refinement}, we present the results of a temporal refinement study for various values of $\epsilon$. We use a fixed mesh with size $128 \times 128$ for this study, and compare the solutions against a reference solution obtained using the same mesh and a CFL number of $0.001$. The methods exhibit third-order accuracy in time, with convergence saturation occurring at smaller CFL numbers, where the spatial error becomes dominant. During this study, we also recorded the average CUR and SVD ranks for each case. As shown in \Cref{fig:consistent temporal refinement rank data}, the average solution rank remains largely unchanged across the different CFL numbers. However, we observe that smaller values of $\epsilon$ tend to yield smaller average ranks, as the system transitions to equilibrium more quickly. The average number of Newton and Krylov iterations at each stage of the time integration are plotted as functions of the CFL in \Cref{fig:consistent temporal refinement iteration data}. The results indicate that smaller values of $\epsilon$ generally require slightly more iterations, but the increase is not substantial. However, as the CFL number increases, the number of nonlinear and linear iterations tends to rise. For moderate CFL numbers (e.g., CFL $< 10$), the iteration counts remain small, so preconditioning for the inner Krylov method may not be necessary. For larger CFL numbers, preconditioning will be essential to ensure the robustness of the solver. This aspect of the methods will be explored in future work.

In \Cref{tab:consistent space refinement}, we present the results of a spatial refinement study in which the mesh resolution is successively doubled from $N_x = 16$ to $N_x = 128$, at a fixed CFL number of $4$, for several values of the Knudsen number $\epsilon$. Third-order convergence is observed in the $L^1$ error for both the kinetic and fluid regimes, which is consistent with the third-order DIRK scheme employed in the proposed method. Achieving fifth-order spatial accuracy would require a smaller CFL number together with tighter truncation tolerances, which would substantially increase the storage cost. Columns five through seven report the absolute errors in the conserved quantities (density $\rho$, momentum $\rho u$, and energy $E$), all of which are preserved to machine precision. Columns eight and nine give the average numerical rank from the SVD and CUR approximations, which remains modest across all cases. The final two columns show the average number of Newton iterations and Krylov iterations per Runge–Kutta stage, both of which remain small.

\begin{figure}[!htbp]
    \centering
    \includegraphics[width=0.325\linewidth]{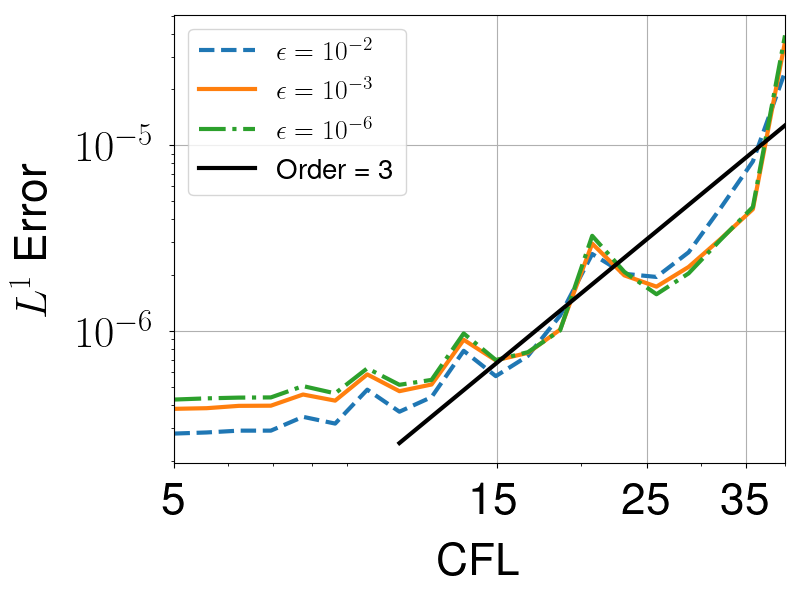}
    \caption{(Consistent initial data) Convergence results for a temporal refinement study.}
    \label{fig:consistent temporal refinement}
\end{figure}

\begin{figure}[!htbp]
    \centering
    \begin{subfigure}[b]{0.325\textwidth}
        \centering
        \includegraphics[width=\textwidth]{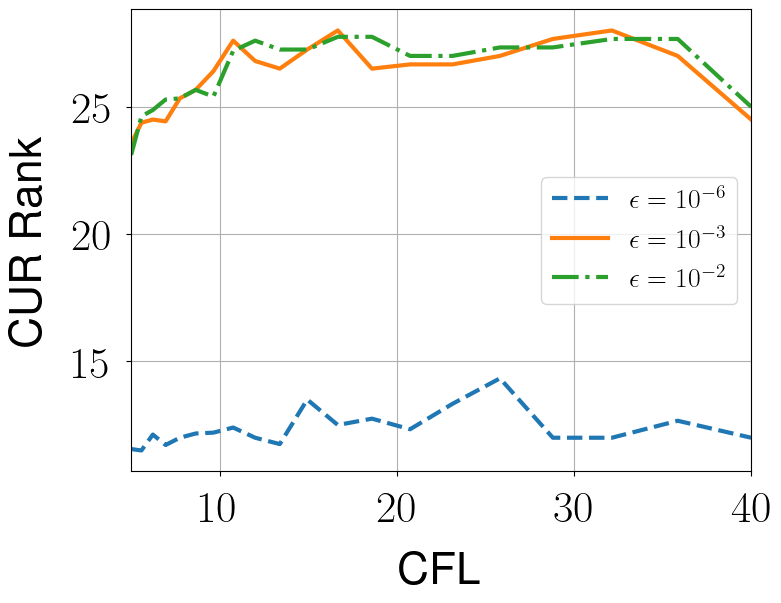}
    \end{subfigure}
    \hspace{1.5em}
    \begin{subfigure}[b]{0.325\textwidth}
        \centering
        \includegraphics[width=\textwidth]{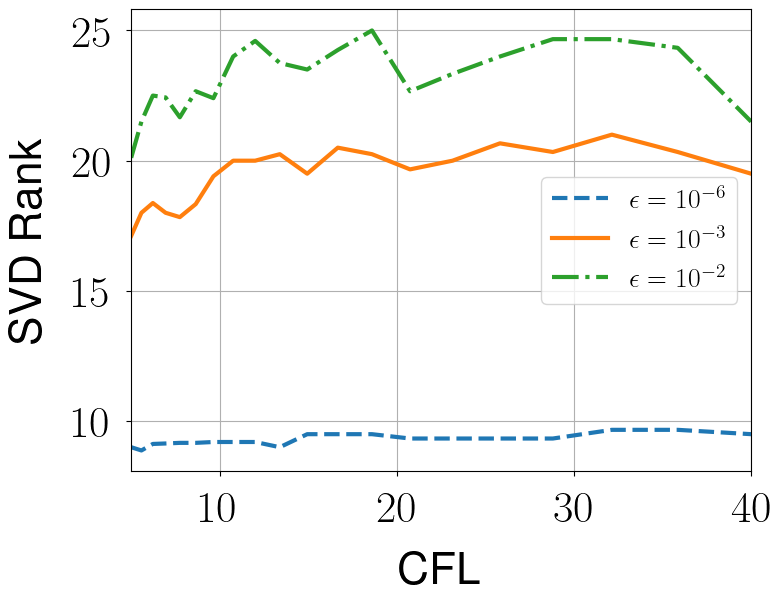}
    \end{subfigure}
    \caption{(Consistent initial data). CUR (left) and SVD (right) rank versus CFL for different $\epsilon$.}
    \label{fig:consistent temporal refinement rank data}
\end{figure}

\begin{figure}[!htbp]
\centering
  \begin{subfigure}[b]{0.325\textwidth}
    \centering
    \includegraphics[width=\textwidth]{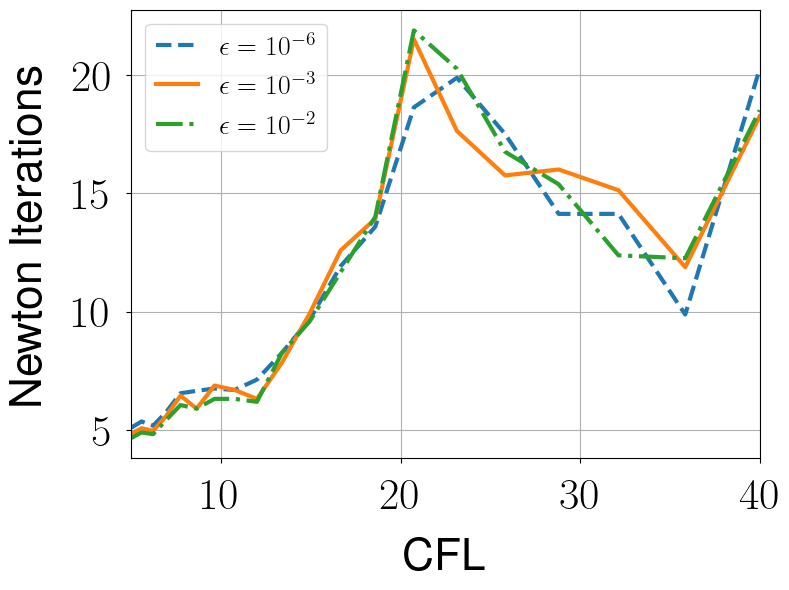}
    \label{fig:newton}
  \end{subfigure}
  \hspace{1.5em}
  \begin{subfigure}[b]{0.325\textwidth}
    \centering
    \includegraphics[width=\textwidth]{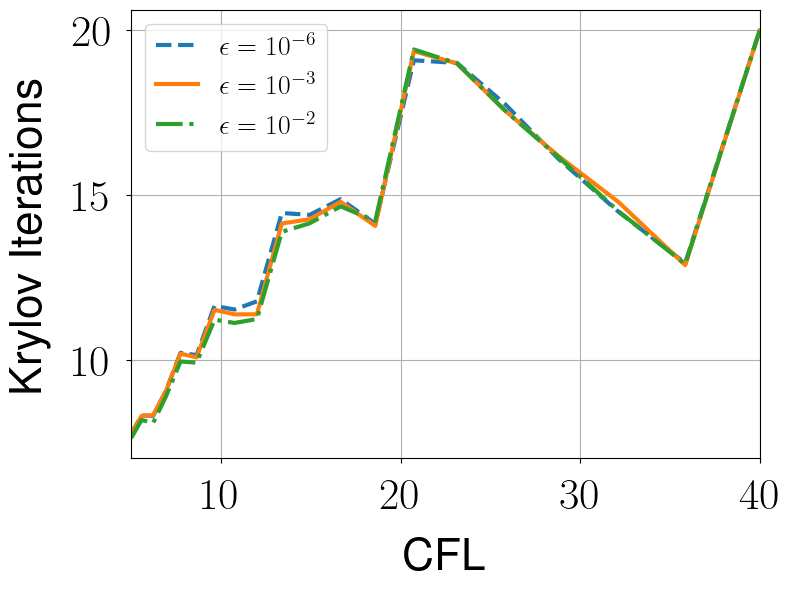}
    \label{fig:krylov}
  \end{subfigure}

  \caption{(Consistent initial data). Average number of Newton (left) and Krylov (right) iterations per stage as a function of the CFL for different $\epsilon$.}
  \label{fig:consistent temporal refinement iteration data}
\end{figure}

\begin{table}[!htbp]
\centering
\caption{(Consistent initial data).  Spatial refinement study using different Knudsen numbers $\epsilon$. The second column reports the number of spatial cells $N_x$. The third and fourth columns provide the $L^1$ approximation error and corresponding convergence order. Columns five through seven list the absolute errors in the conserved quantities (density $\rho$, momentum $\rho u$, and energy $E$). Columns eight and nine show the average numerical rank obtained using SVD and CUR approximations. The final two columns report the average number of outer Newton iterations and inner Krylov iterations per Runge–Kutta stage.}\label{tab:consistent space refinement}
\centering
\begin{tabular}{|c|c||cc||c|c|c||c|c|c|c|}
\hline
$\epsilon$ & $N_{x}$ & $L^1$ error & Order & $\rho$ & $\rho u$ & $E$ & SVD & CUR & Newton & Krylov \\
\hline
&$16$ & 3.0e-04 & --- & 2e-16 & 3e-18 & 2e-16 & 6.0 & 7.0& 4.4& 5.5 \\
$10^{-2}$&$32$ & 4.2e-05 & 2.84 & 4e-16 & 0e+00 & 0e+00 & 9.5 & 11.2& 3.8& 6.1 \\
&$64$& 4.8e-06 & 3.14 & 4e-16 & 3e-18 & 2e-16 & 15.2 & 18.0& 3.5& 6.9 \\
&$128$& 2.8e-07 & 4.10 & 4e-16 & 0e+00 & 2e-16 & 20.4 & 22.8& 4.1& 6.8 \\
\hline
&$16$ & 6.6e-04 & --- & 2e-16 & 3e-18 & 0e+00 & 5.7 & 7.0& 4.6& 5.5 \\
$10^{-3}$&$32$ & 4.9e-05 & 3.73 & 4e-16 & 3e-18 & 0e+00 & 9.2 & 11.0& 4.6& 5.5 \\
&$64$& 5.7e-06 & 3.11 & 4e-16 & 3e-18 & 2e-16 & 14.5 & 18.3& 3.9& 7.1 \\
&$128$& 3.8e-07 & 3.92 & 4e-16 & 3e-18 & 2e-16 & 16.6 & 22.5& 4.4& 6.9 \\
\hline
&$16$ & 9.1e-04 & --- & 2e-16 & 3e-18 & 0e+00 & 5.7 & 7.0& 5.0& 5.3 \\
$10^{-6}$&$32$ & 5.1e-05 & 4.15 & 4e-16 & 0e+00 & 2e-16 & 8.0 & 10.2& 3.7& 6.2 \\
&$64$& 6.3e-06 & 3.03 & 4e-16 & 0e+00 & 0e+00 & 8.7 & 12.3& 4.0& 7.2 \\
&$128$& 4.2e-07 & 3.88 & 4e-16 & 0e+00 & 2e-16 & 8.9 & 11.6& 4.8& 6.9 \\
\hline
\end{tabular}
\end{table}

\end{exa}

\begin{exa}(Riemann problem \cite{pieraccini2007implicit}) Next, we investigate the capabilities of the proposed methods in capturing smooth rarefaction waves as well as non-smooth solution structures, including contact discontinuities and shocks, found in compressible gas dynamics problems. Consider the following initial distribution, which exhibits discontinuous behavior:
\begin{equation}
    f(x,v,0)=\begin{cases}
        \dfrac{\rho_L}{\sqrt{2\pi T_L}}\exp\left(-\dfrac{(v-u_L)^2}{2T_L}\right)\quad x\in[0,0.5],\\
        \dfrac{\rho_R}{\sqrt{2\pi T_R}}\exp\left(-\dfrac{(v-u_R)^2}{2T_R}\right)\quad x\in[0.5,1],
    \end{cases}
\end{equation}
where $(\rho_L,u_L,T_L) = (2.25,0,1.125)$ and $(\rho_R,u_R,T_R) = (3/7,0,1/6)$. Fixed boundary conditions are applied, using values from the initial condition. The phase space is discretized with a mesh containing $256 \times 256$ cells in space and velocity, respectively. The simulation runs until a final time of $t = 0.16$, with a CFL number of $4$. The relative truncation tolerances for the CUR and SVD approximations are set to $\epsilon_C = 10^{-4}$ and $\epsilon_S = 10^{-3}$.

The macroscopic density, velocity, and temperature at the final time step are presented in \Cref{fig:Riemann moments}. The lines in each plot correspond to different Knudsen numbers, specifically $\epsilon = 10^{-6}$, $\epsilon = 10^{-3}$, and $\epsilon = 10^{-2}$. Our results align well with those from \cite{li2023high}, with sharper solution features emerging as $\epsilon$ decreases. These non-smooth features are consistent with the behavior expected from the compressible Euler system. 

\begin{figure}[!htbp]
\centering
  \begin{subfigure}[b]{0.32\textwidth}
    \centering
    \includegraphics[width=\textwidth]{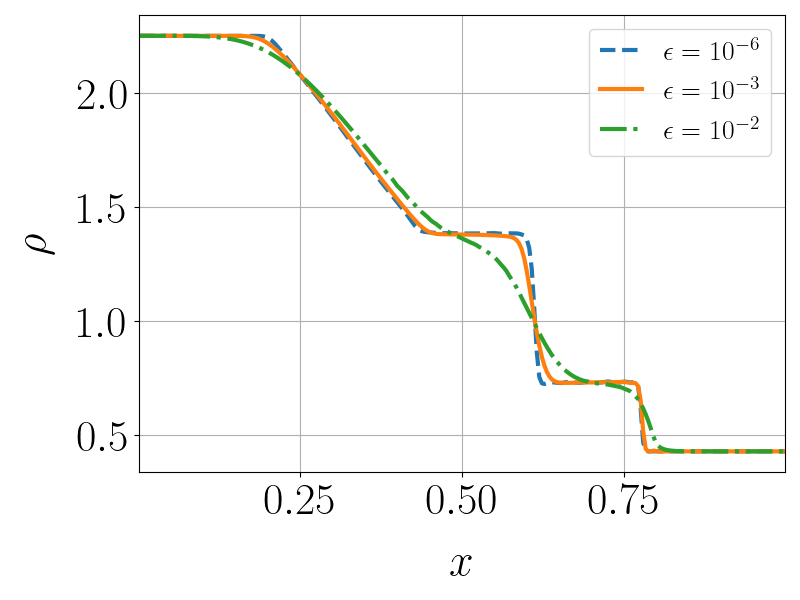}
    \label{fig:Riemann_rho}
  \end{subfigure}
  \begin{subfigure}[b]{0.32\textwidth}
    \centering
    \includegraphics[width=\textwidth]{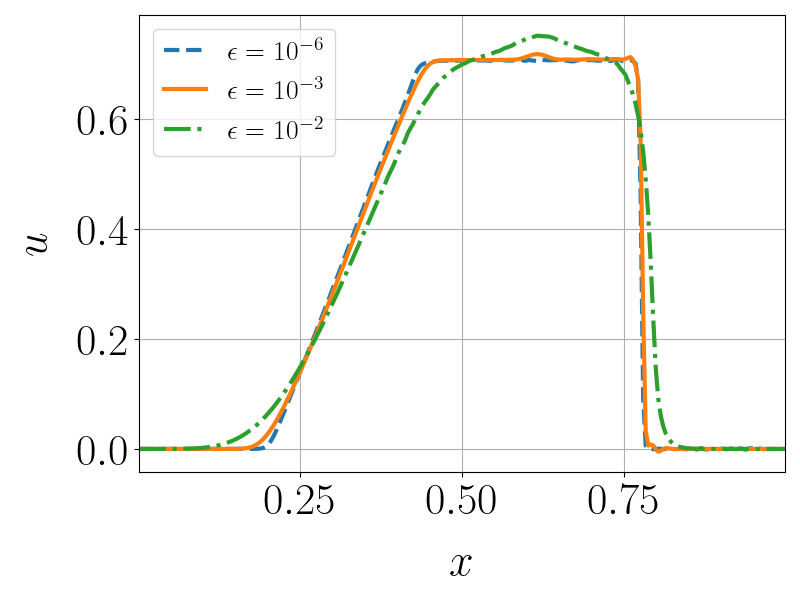}
    \label{fig:Riemann_u}
  \end{subfigure}
  \begin{subfigure}[b]{0.32\textwidth}
    \centering
    \includegraphics[width=\textwidth]{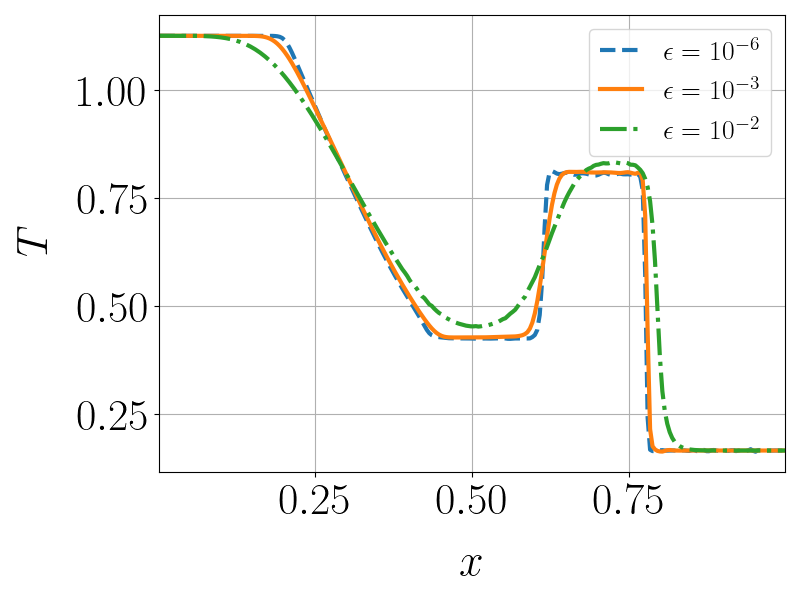}
    \label{fig:Riemann_T}
  \end{subfigure}
  \caption{(Riemann problem). The macroscopic quantities $\rho$ (left), $u$ (middle), and $T$ (right) at time $t = 0.16$ using different values of $\epsilon$.}
  \label{fig:Riemann moments}
\end{figure}

\Cref{fig:Riemann rank data} presents the SVD and CUR rank data as functions of time. As $\epsilon$ decreases, the rank of the solution decreases as well, with this effect being more pronounced in the SVD data. For $\epsilon = 10^{-2}$, the SVD and CUR ranks initially increase at earlier times before reaching a plateau. This increase likely corresponds to the rapid evolution of solution features that occur shortly after the interaction of the Riemann states.

\begin{figure}[!htbp]
    \centering
    \begin{subfigure}[b]{0.325\textwidth}
        \centering
        \includegraphics[width=\textwidth]{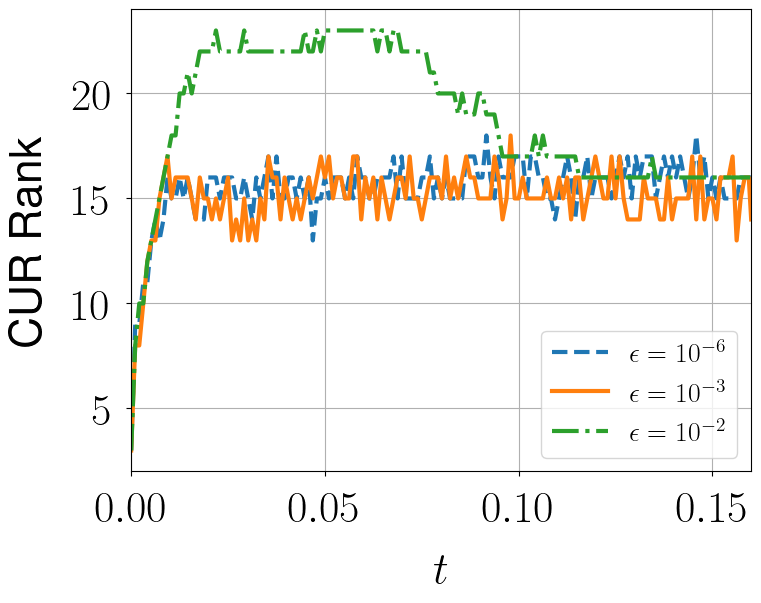}
    \end{subfigure}
    \hspace{1.5em}
    \begin{subfigure}[b]{0.325\textwidth}
        \centering
        \includegraphics[width=\textwidth]{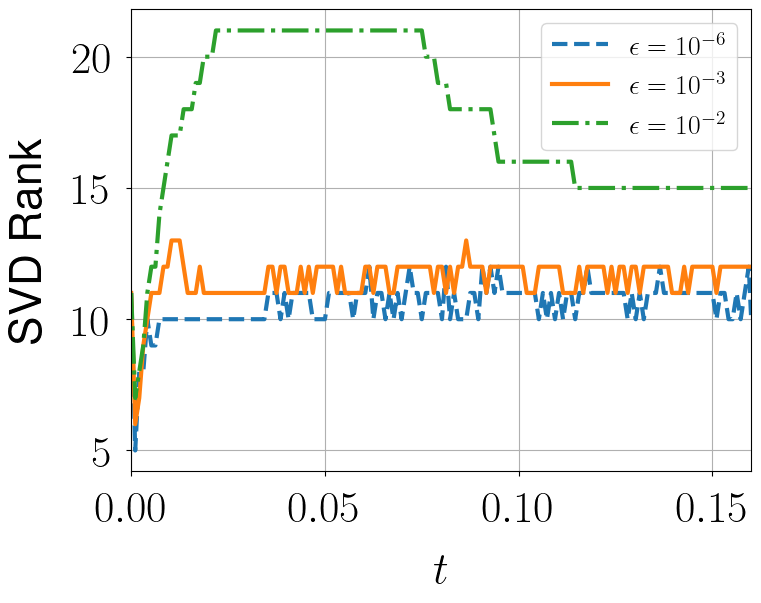}
    \end{subfigure}
    \caption{(Riemann problem). CUR (left) and SVD (right) solution rank versus time for different values of $\epsilon$.}
    \label{fig:Riemann rank data}
\end{figure}

The average number of Newton and Krylov iterations at each stage of the time integration is shown as a function of time in \Cref{fig:Riemann iteration data}. Our results indicate that smaller values of $\epsilon$ generally lead to an SLightly higher number of iterations, but the increase is not substantial. Again, the iteration counts remain small across all cases, even though no preconditioning strategy is applied to the inner Krylov method.

\begin{figure}[!htbp]
    \centering
    \begin{subfigure}[b]{0.325\textwidth}
        \centering
        \includegraphics[width=\textwidth]{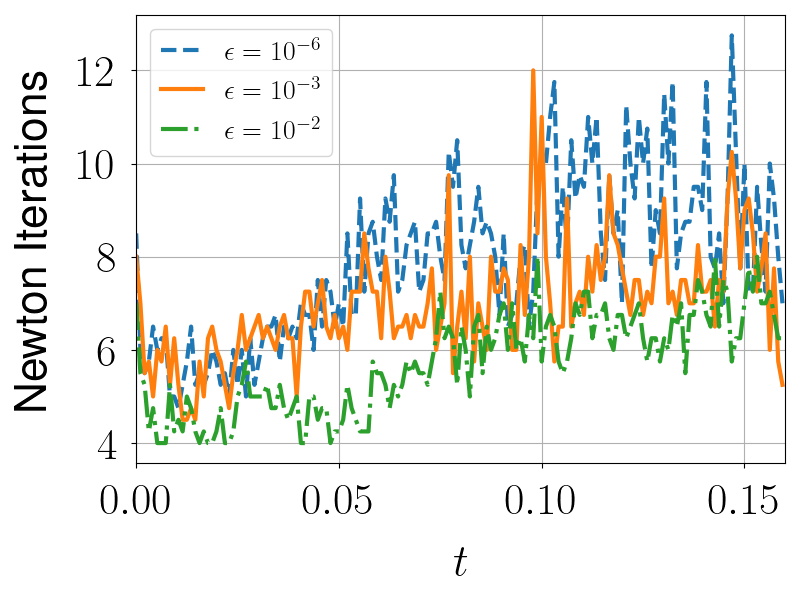}
    \end{subfigure}
    \hspace{1.5em}
    \begin{subfigure}[b]{0.325\textwidth}
        \centering
        \includegraphics[width=\textwidth]{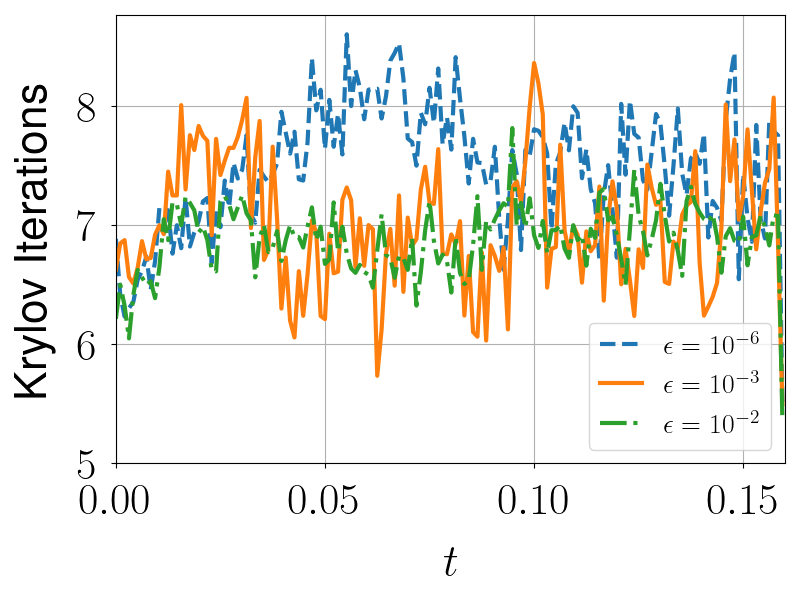}
    \end{subfigure}
    \caption{(Riemann problem). Average number of Newton (left) and Krylov (right) iterations per stage as a function of time obtained with different values of $\epsilon$.}
    \label{fig:Riemann iteration data}
\end{figure}

\end{exa}

\begin{exa}(Mixed regime problem \cite{xiong2015high}) We now apply the proposed method to a more challenging mixed-regime problem to evaluate its ability to handle multi-scale features. We consider a spatially dependent Knudsen number, denoted as $\epsilon(x)$, defined as
\begin{equation}
    \label{eq:variable knudsen number function}
    \epsilon(x) = \epsilon_{0} + \frac{1}{2}\Big(\tanh(1-a_0x)+\tanh(1+a_0x)\Big),
\end{equation}
where $\epsilon_{0} = 10^{-6}$ and $a_{0}$ is a scalar parameter that effectively tunes the width of $\epsilon(x)$. We explore two specific configurations where the transitions between the fluid and kinetic limits are either slow ($a_{0} = 11$) or fast ($a_{0} = 40$). A plot of the Knudsen number defined by \eqref{eq:variable knudsen number function} is shown in \Cref{fig:variable epsilon plot} with these two values of $a_0$. The limiting regimes are separated by approximately six orders of magnitude, and include a full range of intermediate regimes. While methods can be developed to treat the distinct limits in this problem, the proposed method offers a simplified approach to simultaneously address each of the limits, including the intermediate regimes.

\begin{figure}[!htbp]
    \centering
    \includegraphics[width=0.325\linewidth]{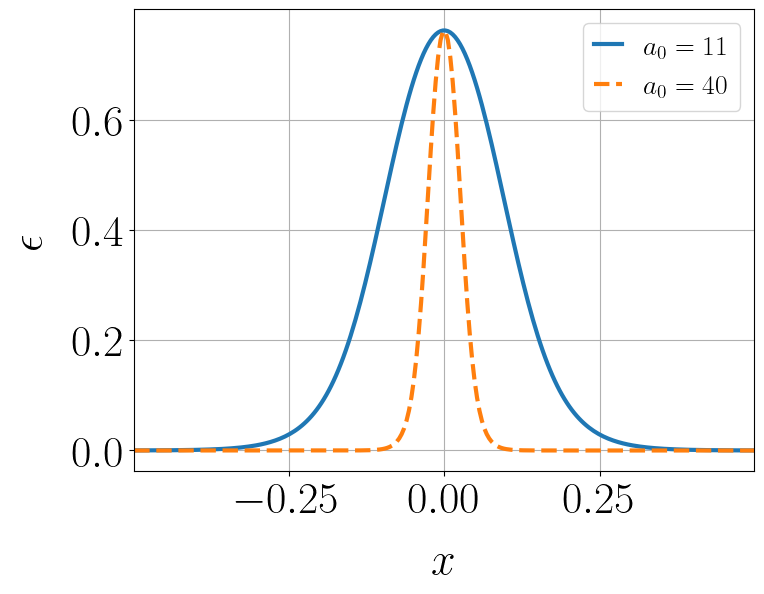}
    \caption{(Mixed regime problem). The Knudsen number $\epsilon(x)$ obtained with $a_{0} = 11$ and $a_{0} = 40$.}
    \label{fig:variable epsilon plot}
\end{figure}

The consistent initial distribution for this problem is given by
\begin{equation*}
    f(x,v,0) = \frac{\rho_{0}}{2\sqrt{2\pi T_{0}}}\left[\exp\left(-\frac{(v-u_{0})^2}{2T_{0}}\right)+\exp\left(-\frac{(v+u_{0})^2}{2T_{0}}\right)\right],\quad x\in[-0.5,0.5],
\end{equation*}
where
\begin{equation*}
    \rho_{0}(x) = 1+0.875\sin(2\pi x), \quad u_{0}(x) = 0.75, \quad T_{0}(x) =0.5+0.4\sin(2\pi x).
\end{equation*}
We apply periodic boundary conditions and discretize phase space using a $256 \times 256$ mesh. The final simulation time is $t = 0.45$, and each configuration is considered for three different CFL numbers: $1$, $1.5$, and $2$. We use $\epsilon_C = 10^{-8}$ and $\epsilon_S = 10^{-7}$ as the relative truncation tolerances for the CUR and SVD approximations.

Plots of the macroscopic variables $\rho$, $u$, and $T$, at $t = 0.1$, $t = 0.3$, and $t = 0.45$, using $a_0 = 11$ and $a_0 = 40$, are shown in \Cref{fig:variable_Knudsen_a0_11_rho_u_T} and \Cref{fig:variable_Knudsen_a0_40_rho_u_T}, respectively, using different CFL numbers. In general, for $a_0 = 11$, the solution exhibits simple, localized structures with some non-smooth features, such as jumps and kinks. In contrast, for $a_0 = 40$, the rapid transition between fluid and kinetic regimes leads to more complex dynamics, in particular, larger jumps and discontinuities. Interestingly, the results across different CFL numbers are visually indistinguishable, even at later times, despite the larger time steps in simulations with higher CFL numbers.

\begin{figure}[!htbp]
\centering
    \begin{subfigure}[b]{0.3\textwidth}
        \centering
        \includegraphics[width=\textwidth,clip,trim={0cm 0cm 0cm 0cm}]{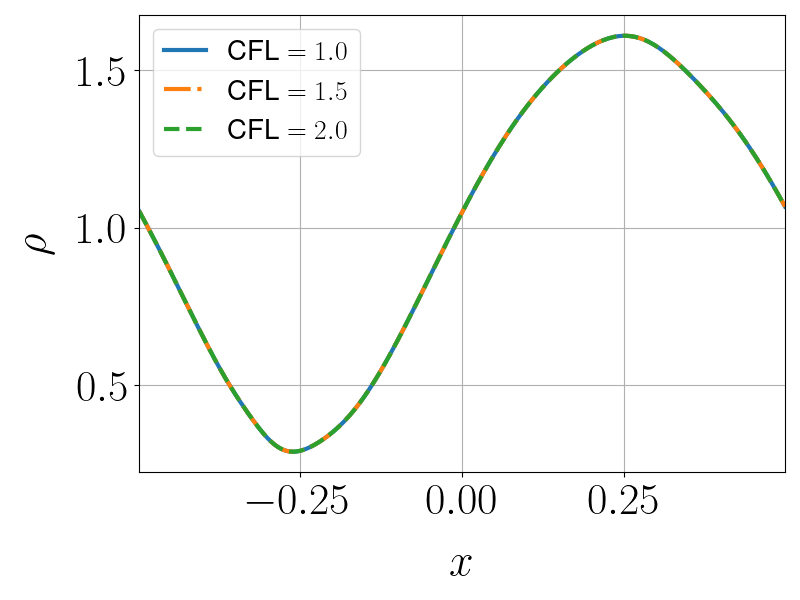}
    \end{subfigure}
    \hspace{-0.5em}
    \begin{subfigure}[b]{0.3\textwidth}
        \centering
        \includegraphics[width=\textwidth,clip,trim={0cm 0cm 0cm 0cm}]{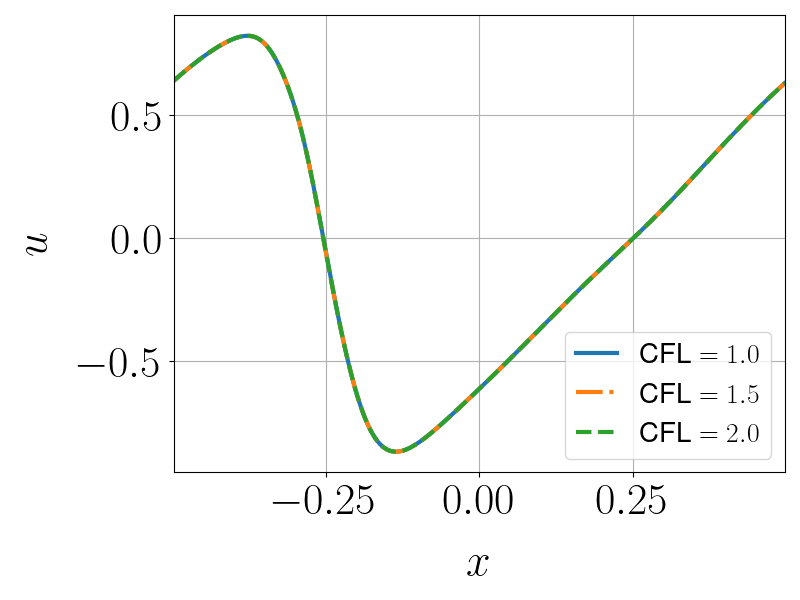}
    \end{subfigure}
    \hspace{-0.5em}
    \begin{subfigure}[b]{0.3\textwidth}
        \centering
        \includegraphics[width=\textwidth,clip,trim={0cm 0cm 0cm 0cm}]{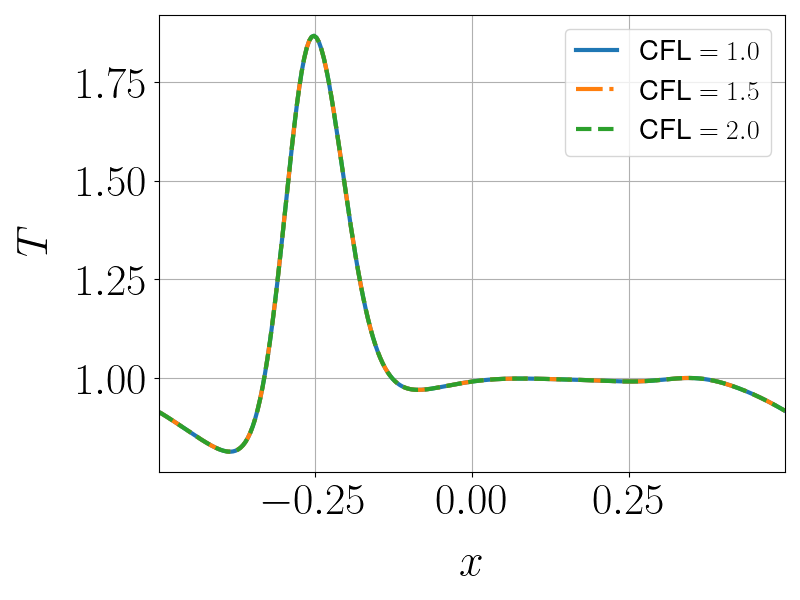}
    \end{subfigure}
    
    \begin{subfigure}[b]{0.3\textwidth}
        \centering
        \includegraphics[width=\textwidth,clip,trim={0cm 0cm 0cm 0cm}]{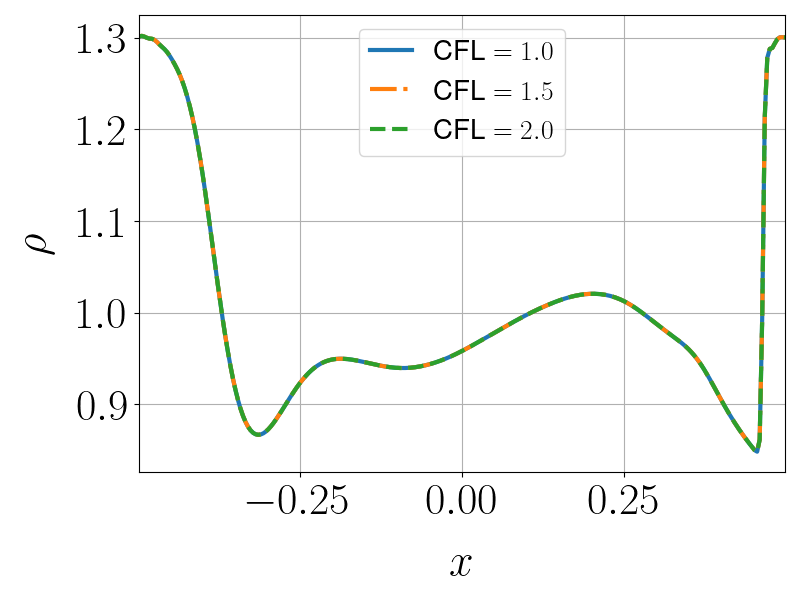}
    \end{subfigure}
    \hspace{-0.5em}
    \begin{subfigure}[b]{0.3\textwidth}
        \centering
        \includegraphics[width=\textwidth,clip,trim={0cm 0cm 0cm 0cm}]{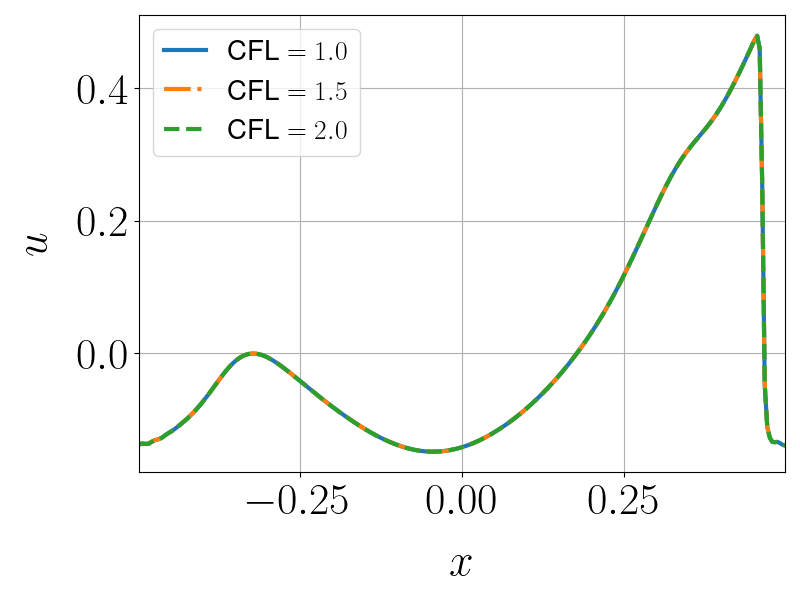}
    \end{subfigure}
    \hspace{-0.5em}
    \begin{subfigure}[b]{0.3\textwidth}
        \centering
        \includegraphics[width=\textwidth,clip,trim={0cm 0cm 0cm 0cm}]{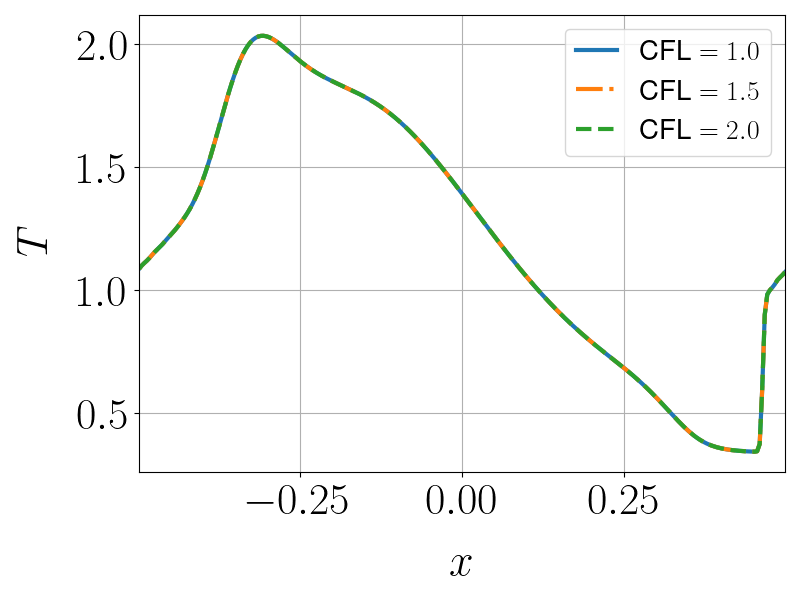}
    \end{subfigure}
    
    \begin{subfigure}[b]{0.3\textwidth}
        \centering
        \includegraphics[width=\textwidth,clip,trim={0cm 0cm 0cm 0cm}]{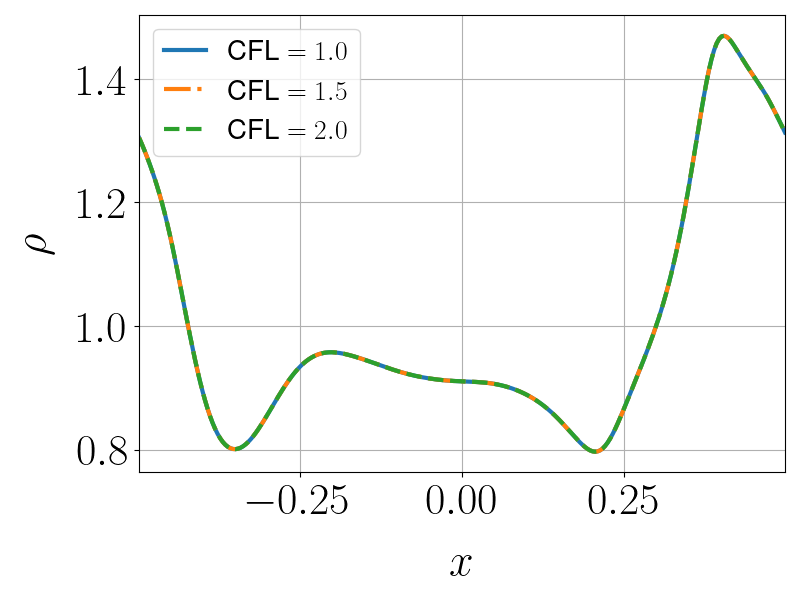}
    \end{subfigure}
    \hspace{-0.5em}
    \begin{subfigure}[b]{0.3\textwidth}
        \centering
        \includegraphics[width=\textwidth,clip,trim={0cm 0cm 0cm 0cm}]{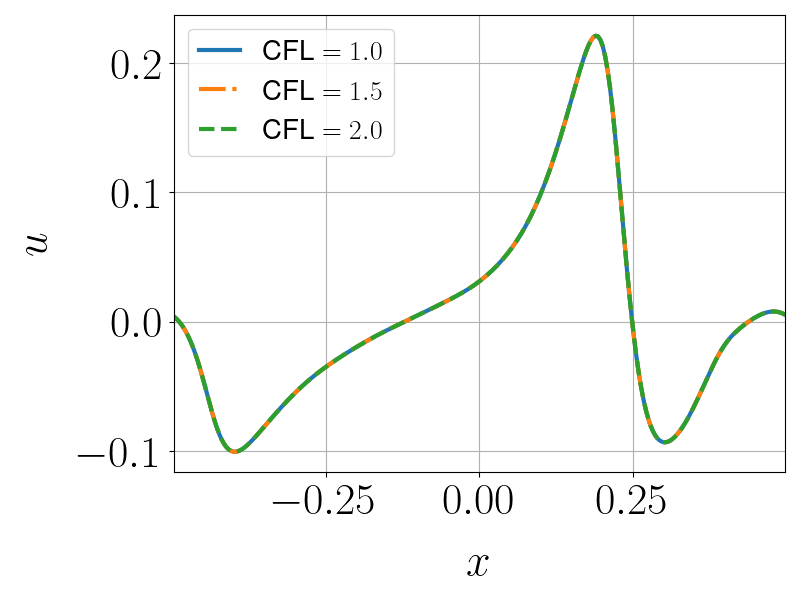}
    \end{subfigure}
    \hspace{-0.5em}
    \begin{subfigure}[b]{0.3\textwidth}
        \centering
        \includegraphics[width=\textwidth,clip,trim={0cm 0cm 0cm 0cm}]{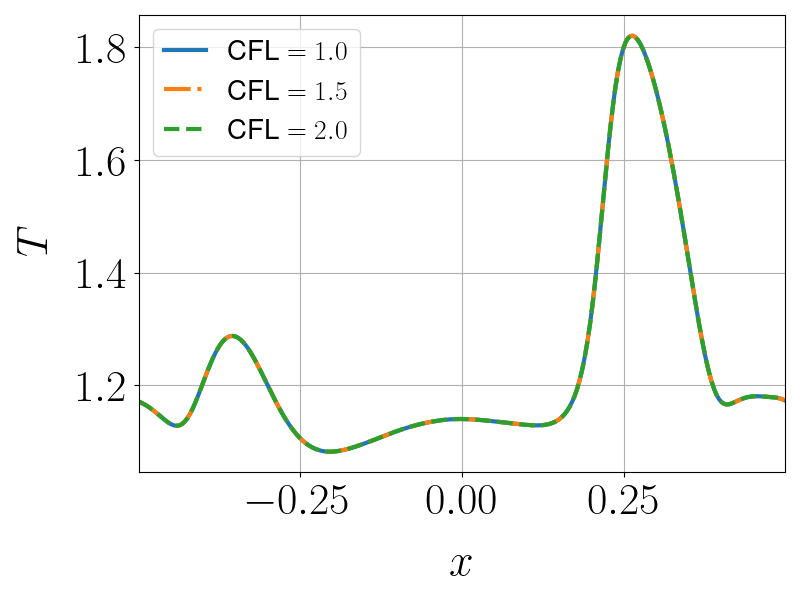}
    \end{subfigure}

    \caption{(Mixed regime problem, $a_{0} = 11$). The macroscopic quantities $\rho$, $u$, and $T$ are shown at time $t = 0.1$ (top row), $t = 0.3$ (middle row), $t=0.45$ (bottom row) with different CFL numbers.}
    \label{fig:variable_Knudsen_a0_11_rho_u_T}
\end{figure}

\begin{figure}[!htbp]
    \centering
    \begin{subfigure}[b]{0.3\textwidth}
        \includegraphics[width=\textwidth,clip,trim={0cm 0cm 0cm 0cm}]{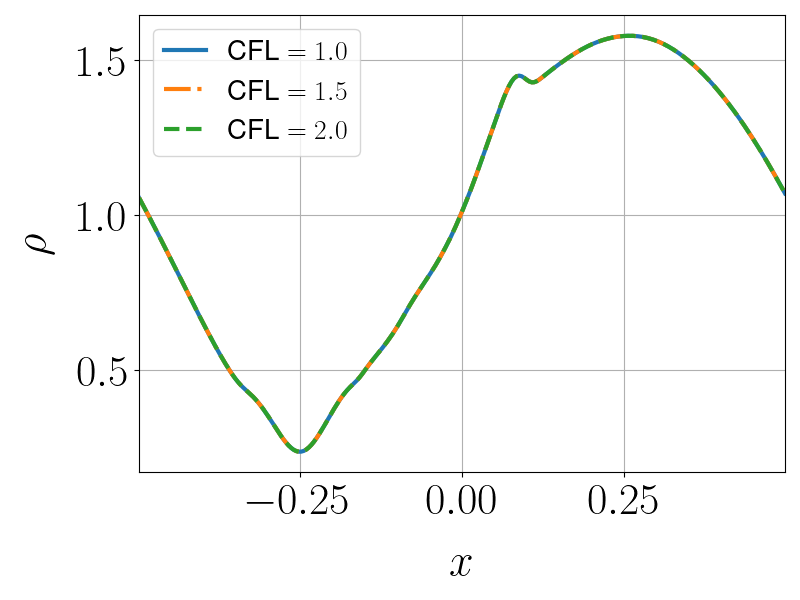}
    \end{subfigure}
    \hspace{-0.5em}
    \begin{subfigure}[b]{0.3\textwidth}
        \includegraphics[width=\textwidth,clip,trim={0cm 0cm 0cm 0cm}]{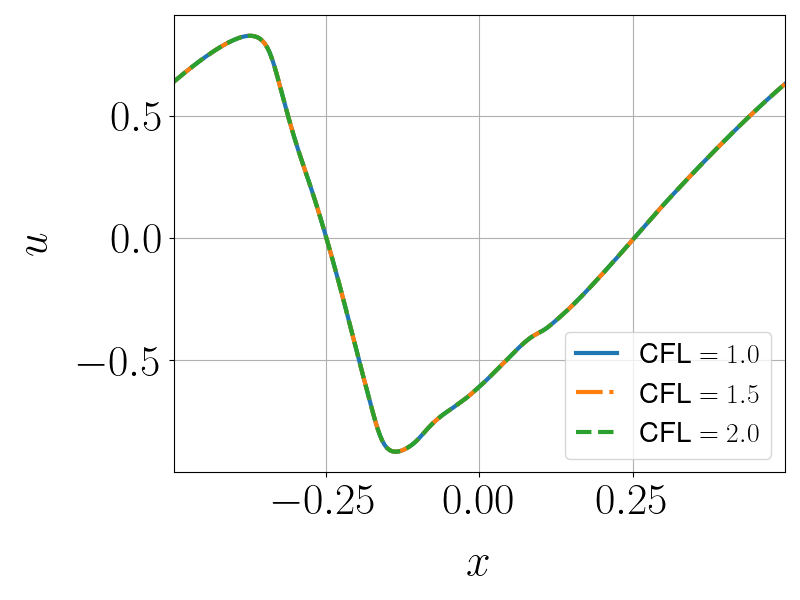}
    \end{subfigure}
    \hspace{-0.5em}
    \begin{subfigure}[b]{0.3\textwidth}
        \includegraphics[width=\textwidth,clip,trim={0cm 0cm 0cm 0cm}]{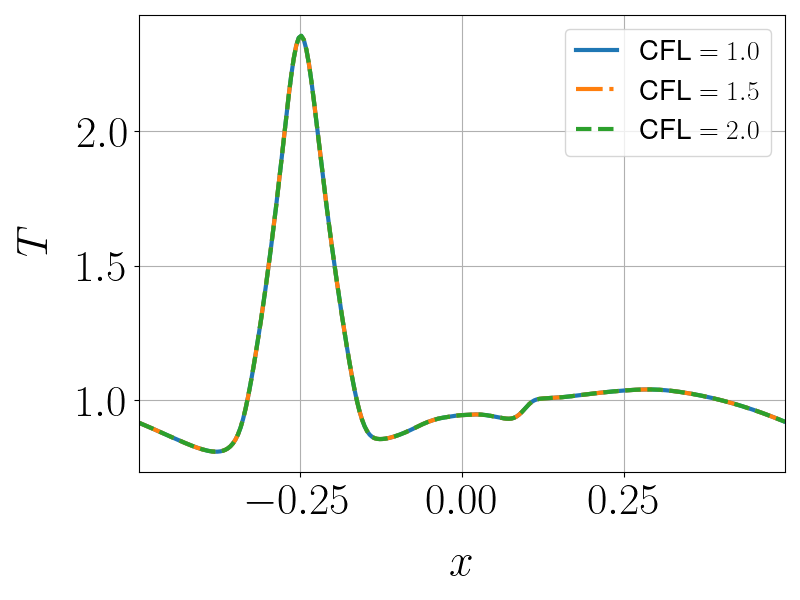}
    \end{subfigure}
    
    \begin{subfigure}[b]{0.3\textwidth}
        \includegraphics[width=\textwidth,clip,trim={0cm 0cm 0cm 0cm}]{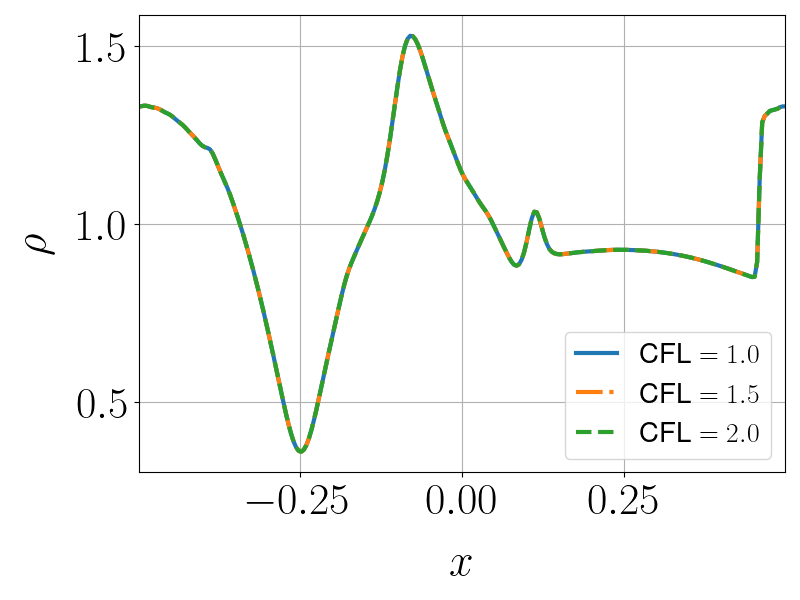}
    \end{subfigure}
    \hspace{-0.5em}
    \begin{subfigure}[b]{0.3\textwidth}
        \includegraphics[width=\textwidth,clip,trim={0cm 0cm 0cm 0cm}]{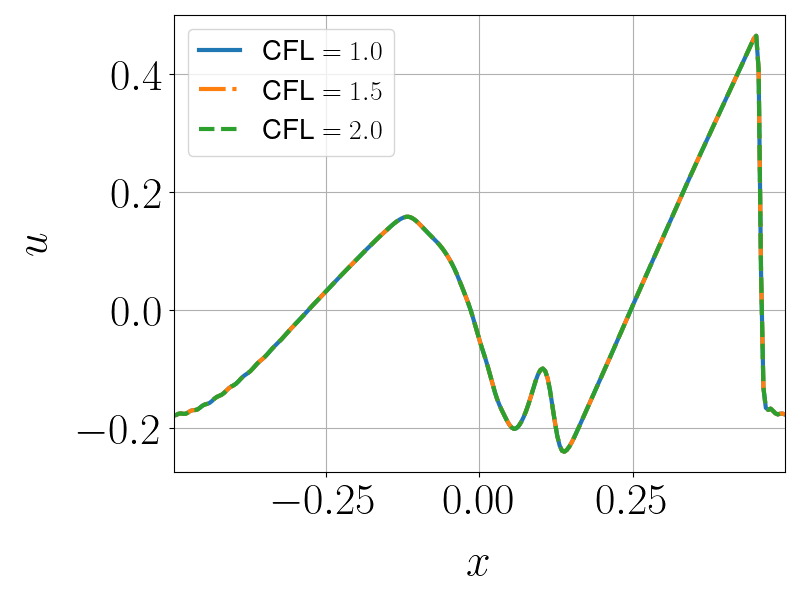}
    \end{subfigure}
    \hspace{-0.5em}
    \begin{subfigure}[b]{0.3\textwidth}
        \includegraphics[width=\textwidth,clip,trim={0cm 0cm 0cm 0cm}]{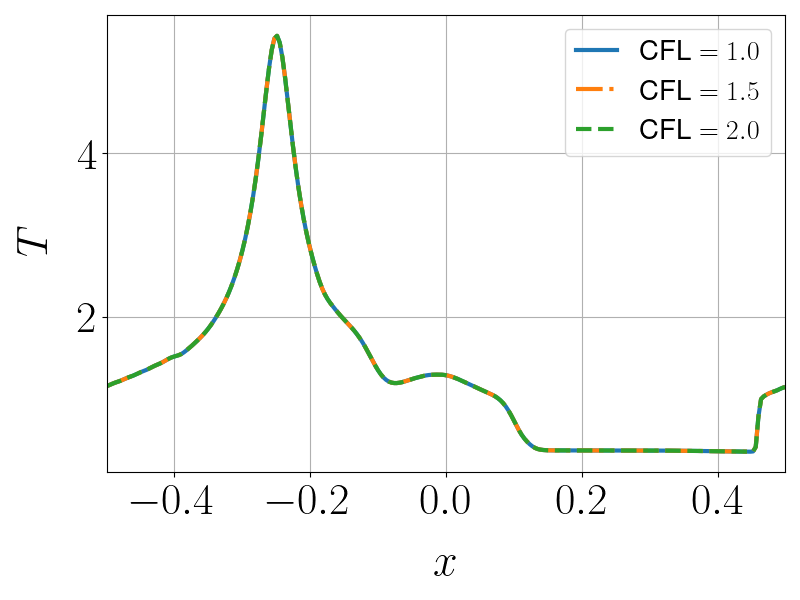}
    \end{subfigure}
    
    \begin{subfigure}[b]{0.3\textwidth}
        \includegraphics[width=\textwidth,clip,trim={0cm 0cm 0cm 0cm}]{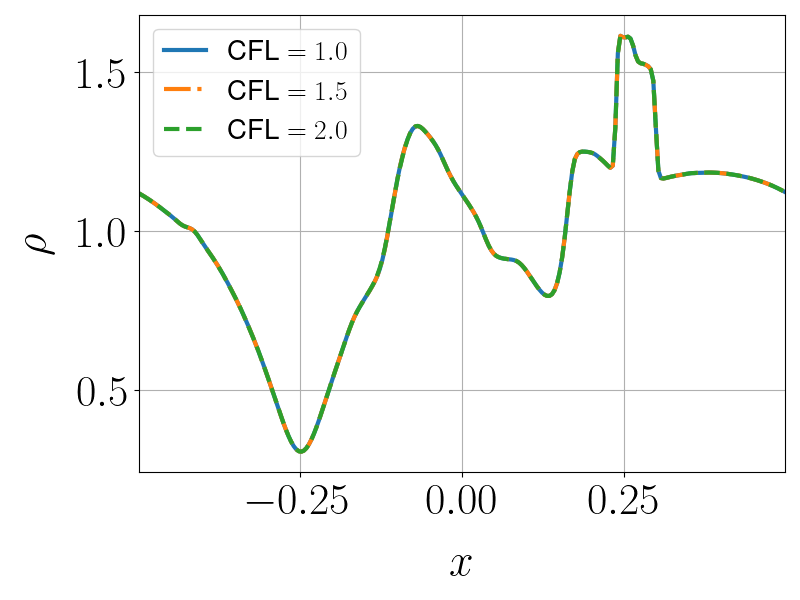}
    \end{subfigure}
    \hspace{-0.5em}
    \begin{subfigure}[b]{0.3\textwidth}
        \includegraphics[width=\textwidth,clip,trim={0cm 0cm 0cm 0cm}]{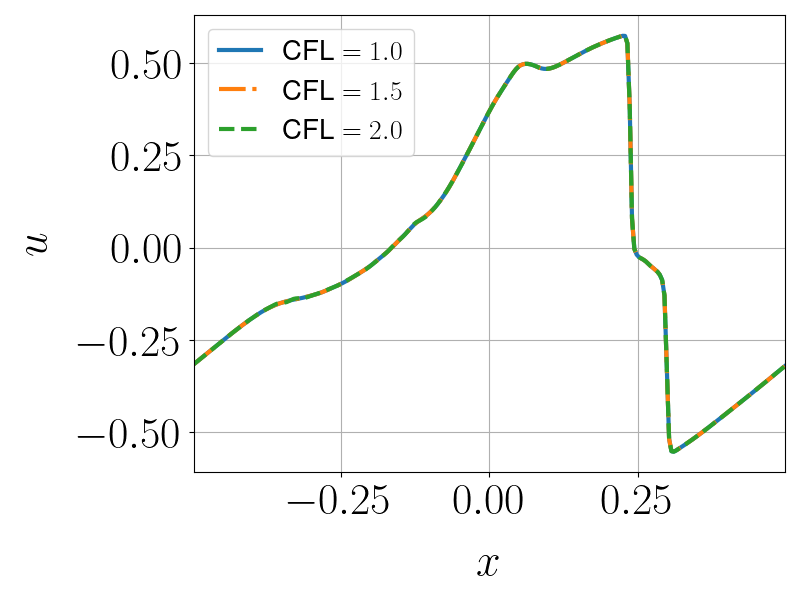}
    \end{subfigure}
    \hspace{-0.5em}
    \begin{subfigure}[b]{0.3\textwidth}
        \includegraphics[width=\textwidth,clip,trim={0cm 0cm 0cm 0cm}]{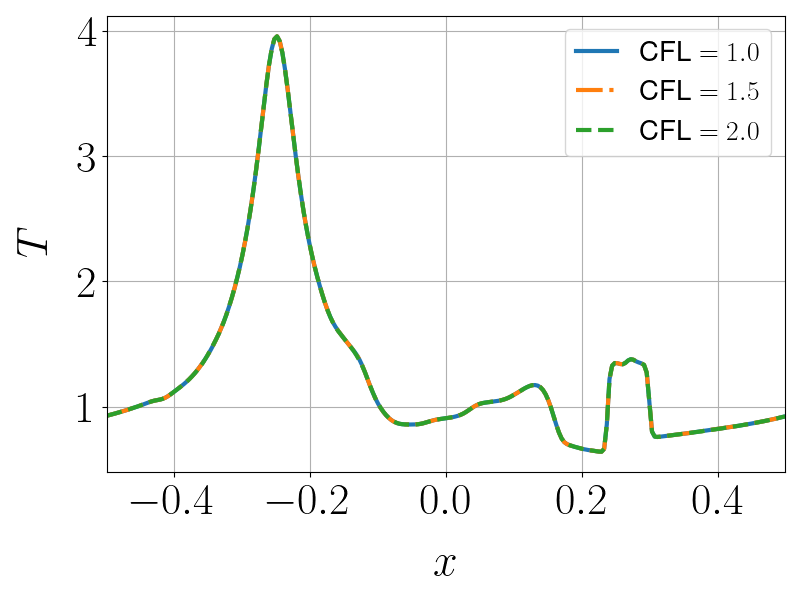}
    \end{subfigure}
    
    \caption{(Mixed regime problem, $a_{0} = 40$). The macroscopic quantities $\rho$, $u$, and $T$ are shown at time $t = 0.1$ (top row), $t = 0.3$ (middle row), $t=0.45$ (bottom row) with different CFL numbers.}
    \label{fig:variable_Knudsen_a0_40_rho_u_T}
\end{figure}

The CUR and SVD ranks for both configurations are presented in \Cref{fig:Variable Knudsen rank data} as a function of time. As expected, the CUR rank is slightly larger than the SVD rank, but both ranks are comparable across the range of CFL values considered. For the case $a_0 = 11$, the solution structures are relatively simple, resulting in small SVD and CUR ranks over time. From this data, we find that the low-rank approach requires approximately 15\% of the storage compared to the full grid method. In contrast, larger ranks are observed for $a_0 = 40$, where more abrupt transitions occur between the fluid and kinetic regimes. Here, the low-rank approach uses approximately 23\% of the storage required by the full grid method, which is higher than the case for $a_0 = 11$. This behavior is expected due to the presence of discontinuities in the solution, as shown in \Cref{fig:variable_Knudsen_a0_40_rho_u_T}. 

\begin{figure}[!htbp]
    \centering
    \begin{subfigure}[b]{0.325\textwidth}
        \includegraphics[width=\textwidth,clip,trim={0cm 0cm 0cm 0cm}]{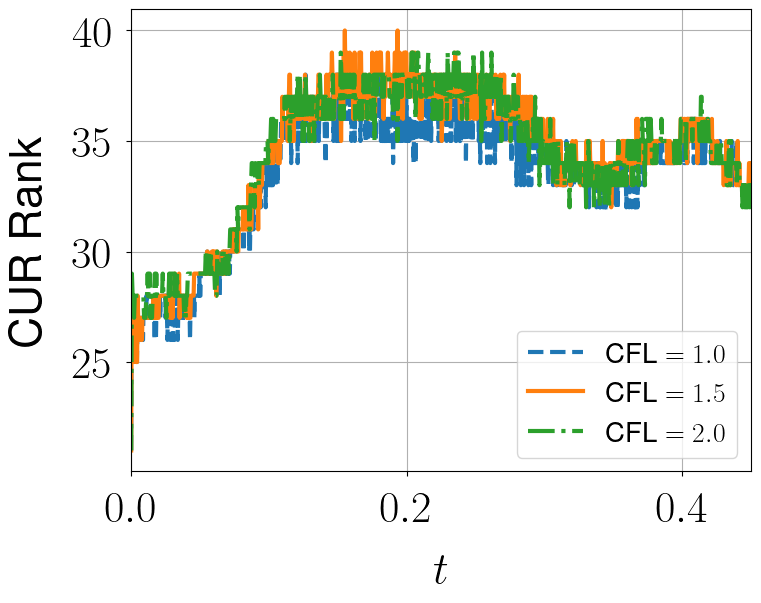}
    \end{subfigure}
    \hspace{0.5em}
    \begin{subfigure}[b]{0.325\textwidth}
        \includegraphics[width=\textwidth,clip,trim={0cm 0cm 0cm 0cm}]{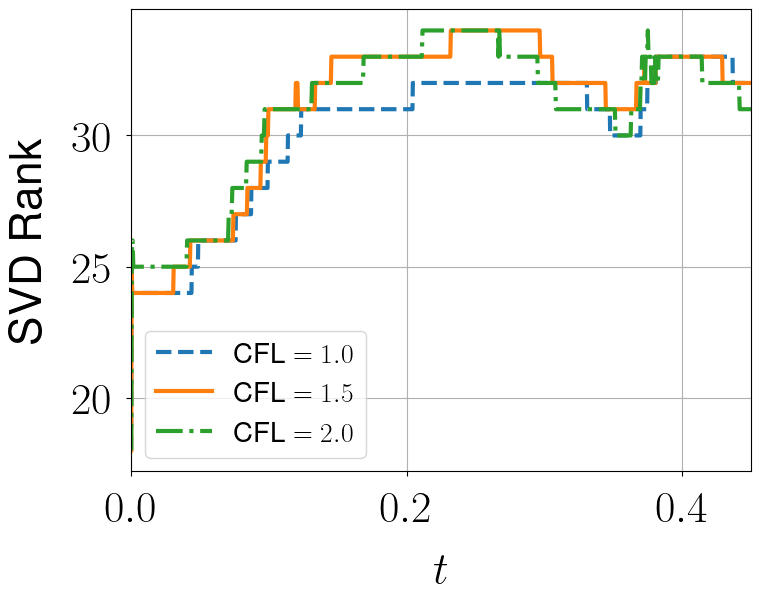}
    \end{subfigure}
    
    \vspace{1em}
    
    \begin{subfigure}[b]{0.325\textwidth}
        \includegraphics[width=\textwidth,clip,trim={0cm 0cm 0cm 0cm}]{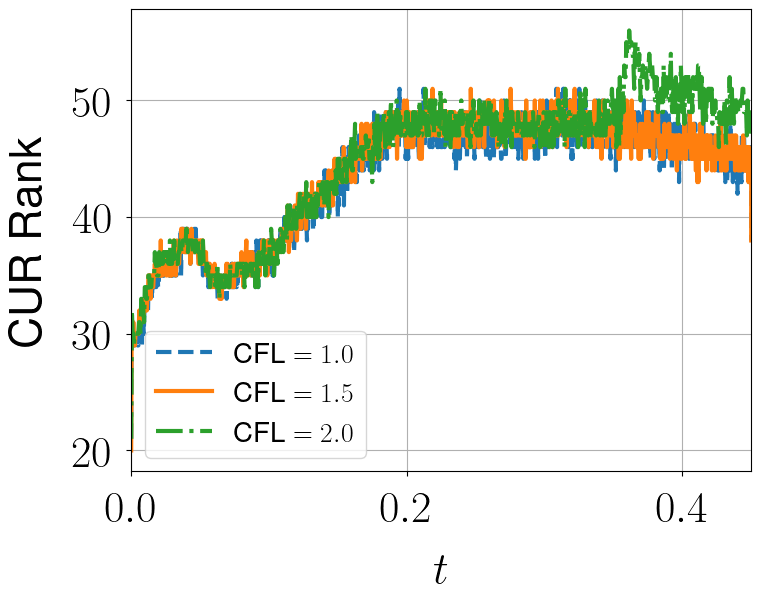}
    \end{subfigure}
    \hspace{0.5em}
    \begin{subfigure}[b]{0.325\textwidth}
        \includegraphics[width=\textwidth,clip,trim={0cm 0cm 0cm 0cm}]{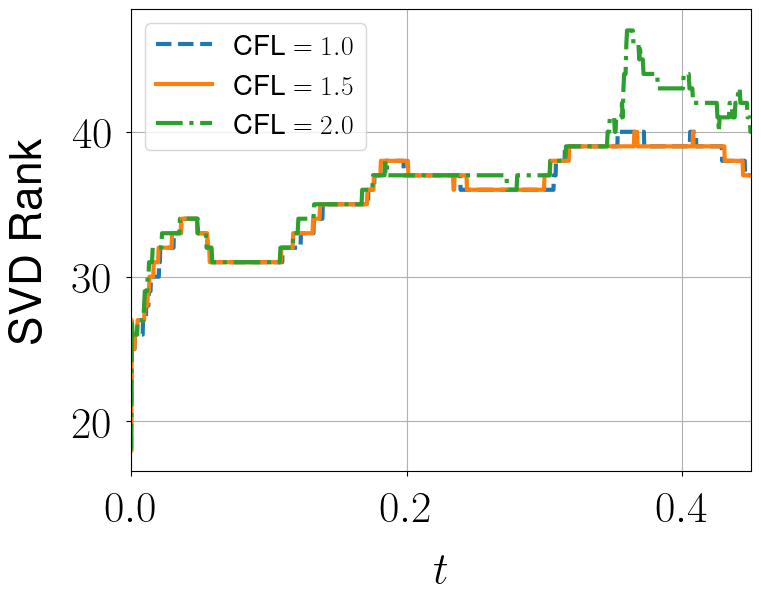}
    \end{subfigure}
    \caption{(Mixed regime problem). CUR and SVD rank over time using different CFL numbers. Figures in the top and bottom rows correspond to $a_{0} = 11$ and $a_{0} = 40$, respectively.}
    \label{fig:Variable Knudsen rank data}
\end{figure}

\Cref{fig:Variable Knudsen iteration data} displays the average number of Newton and Krylov iterations per stage as a function of time. Again, we observe that the number of iterations remains small throughout the simulation, so we find preconditioning to not be necessary. This is primarily due to the high-quality initial guess provided by the low-rank kinetic solver, which significantly improves the convergence of the nonlinear solver. The number of iterations in both cases of $a_{0}$ are generally similar. However, the number of iterations generally increases with the CFL number, as observed in \Cref{fig:consistent temporal refinement iteration data}. In each configuration, the method conserves total mass, momentum, and energy to machine precision, as demonstrated in \Cref{fig:Variable Knudsen conservation data}.

Lastly, \Cref{fig:Variable Knudsen scaling data} presents a plot of the simulation wall time per step versus $N$, the number of mesh points per dimension. A shorter final time of $T = 0.001$ was used with $\text{CFL} = 1$, and the value of $N$ was successively doubles from $N = 8$ to $N = 2048$. In both cases, i.e., $a_{0} = 11$ and $a_{0} = 40$, we observe the correct linear scaling with respect to the mesh size $N$ in the proposed adaptive-rank method.

\begin{figure}[!htbp]
\centering
    \begin{subfigure}[b]{0.325\textwidth}
        \includegraphics[width=\textwidth,clip,trim={0cm 0cm 0cm 0cm}]{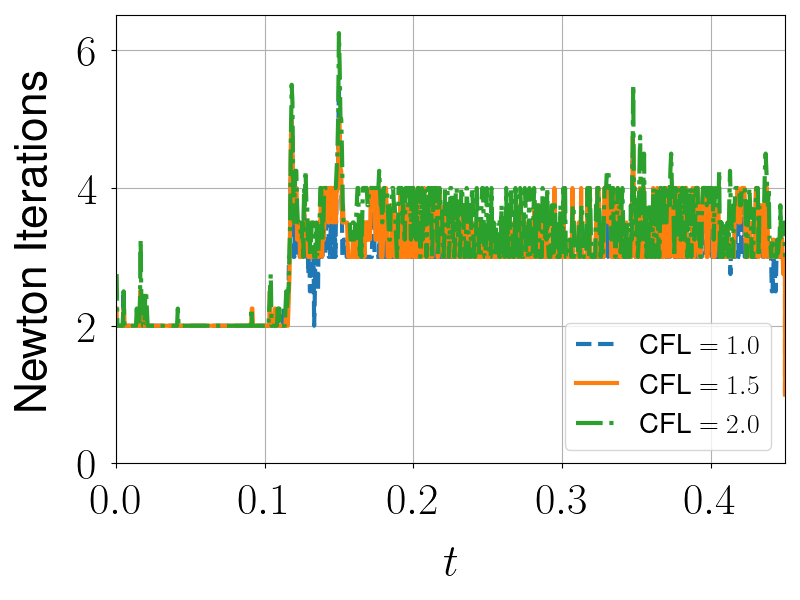}
    \end{subfigure}
    \hspace{0.5em}
    \begin{subfigure}[b]{0.325\textwidth}
        \includegraphics[width=\textwidth,clip,trim={0cm 0cm 0cm 0cm}]{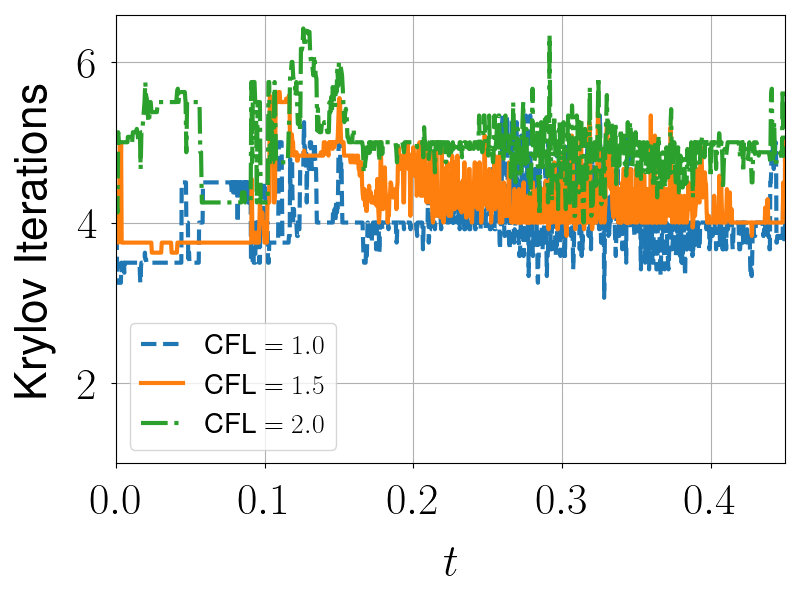}
    \end{subfigure}
    
    \begin{subfigure}[b]{0.325\textwidth}
        \includegraphics[width=\textwidth]{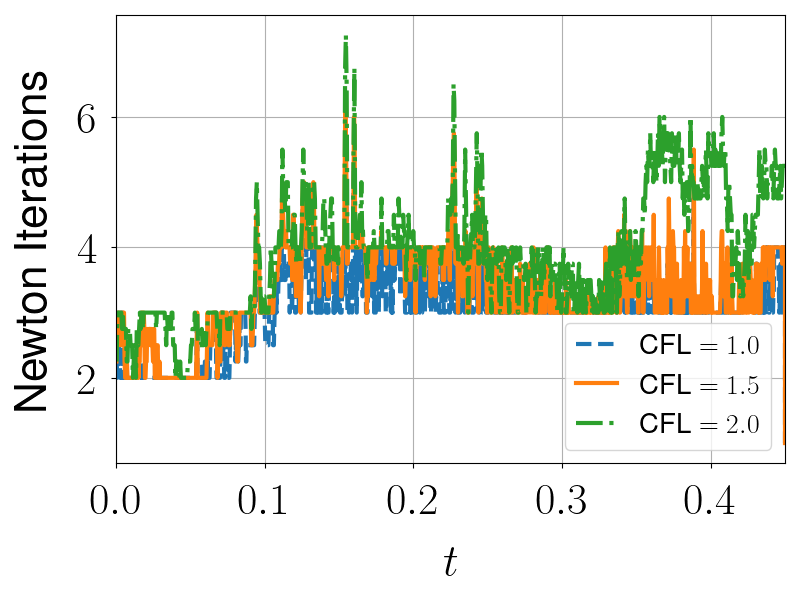}
    \end{subfigure}
    \hspace{0.5em}
    \begin{subfigure}[b]{0.325\textwidth}
        \includegraphics[width=\textwidth]{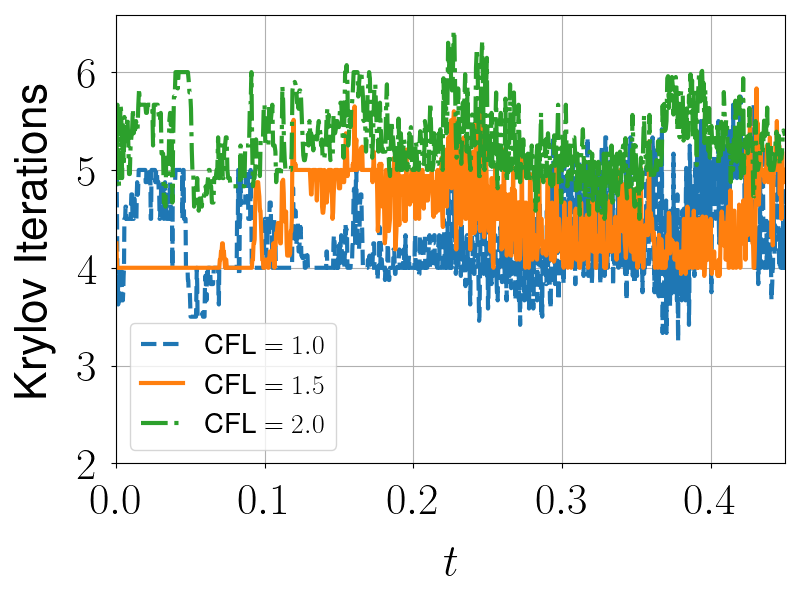}
    \end{subfigure}
    
\caption{(Mixed regime problem). Average number of iterations per stage over time obtained with different CFL numbers. The top and bottom rows represent $a_{0} = 11$ and $a_{0} = 40$, respectively.}
\label{fig:Variable Knudsen iteration data}
\end{figure}

\begin{figure}[!htbp]
\centering
    \begin{subfigure}[b]{0.32\textwidth}
        \includegraphics[width=\textwidth,clip,trim={0cm 0cm 0cm 0cm}]{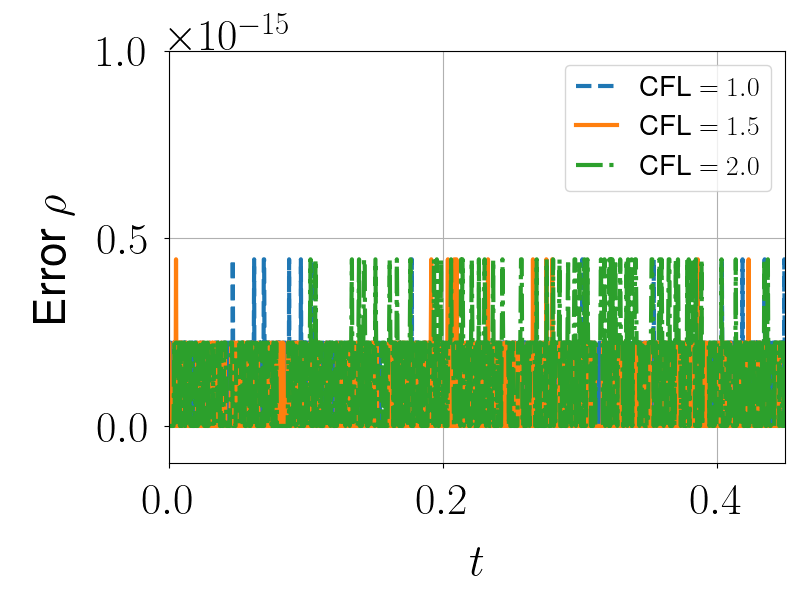}
    \end{subfigure}
    \begin{subfigure}[b]{0.32\textwidth}
        \includegraphics[width=\textwidth,clip,trim={0cm 0cm 0cm 0cm}]{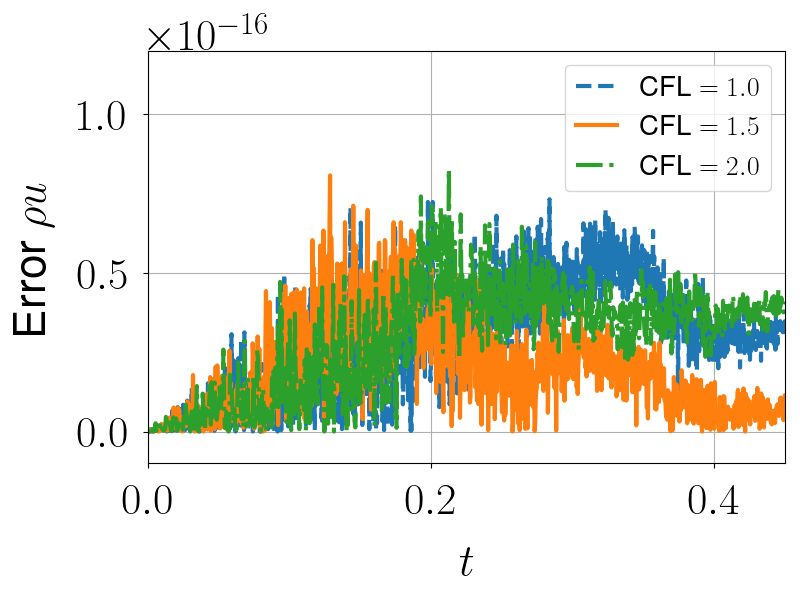}
    \end{subfigure}
    \begin{subfigure}[b]{0.32\textwidth}
        \includegraphics[width=\textwidth,clip,trim={0cm 0cm 0cm 0cm}]{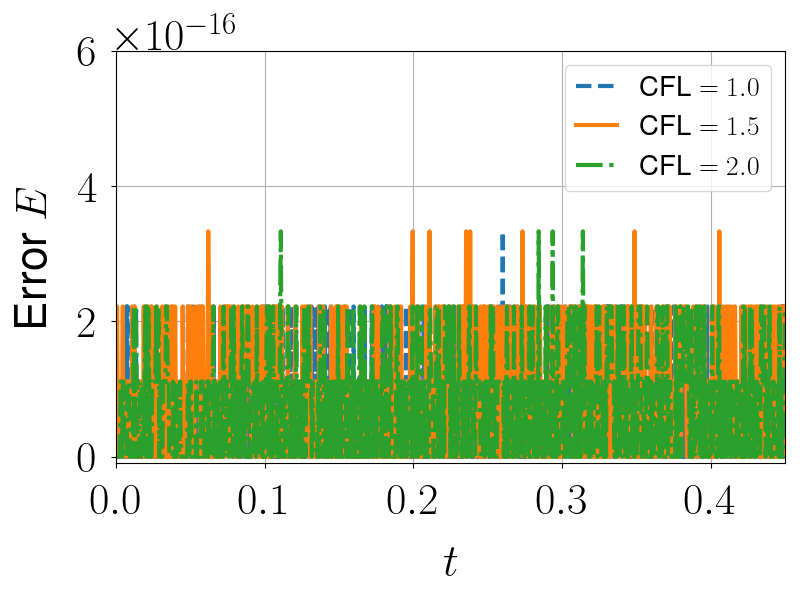}
    \end{subfigure}
    
    \begin{subfigure}[b]{0.32\textwidth}
        \includegraphics[width=\textwidth,clip,trim={0cm 0cm 0cm 0cm}]{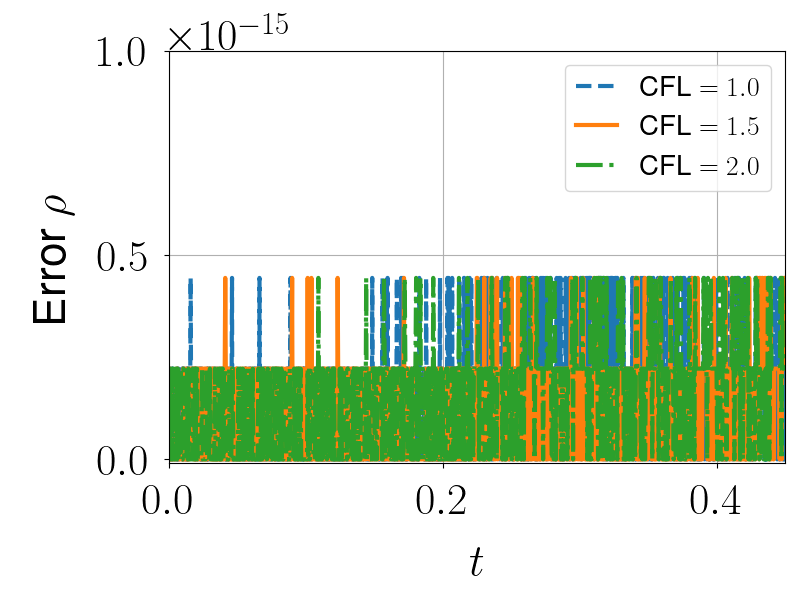}
    \end{subfigure}
    \begin{subfigure}[b]{0.32\textwidth}
        \includegraphics[width=\textwidth,clip,trim={0cm 0cm 0cm 0cm}]{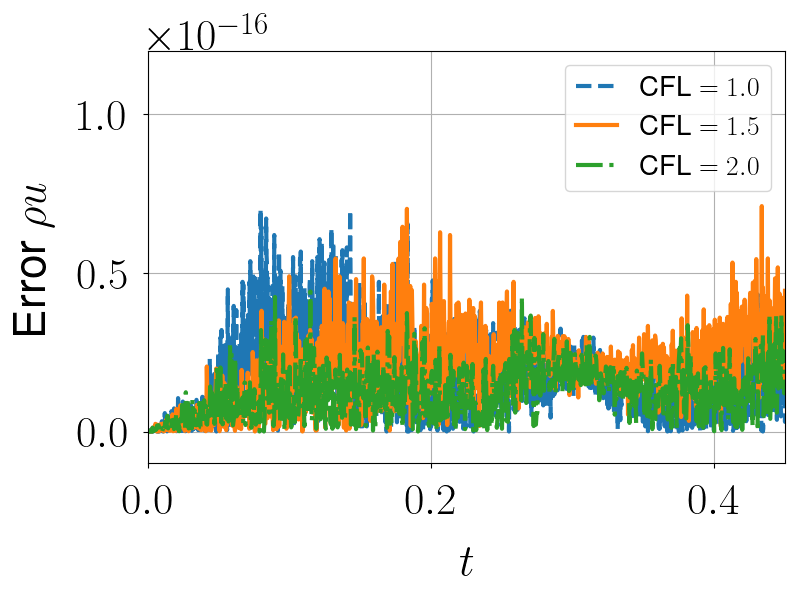}
    \end{subfigure}
    \begin{subfigure}[b]{0.32\textwidth}
        \includegraphics[width=\textwidth,clip,trim={0cm 0cm 0cm 0cm}]{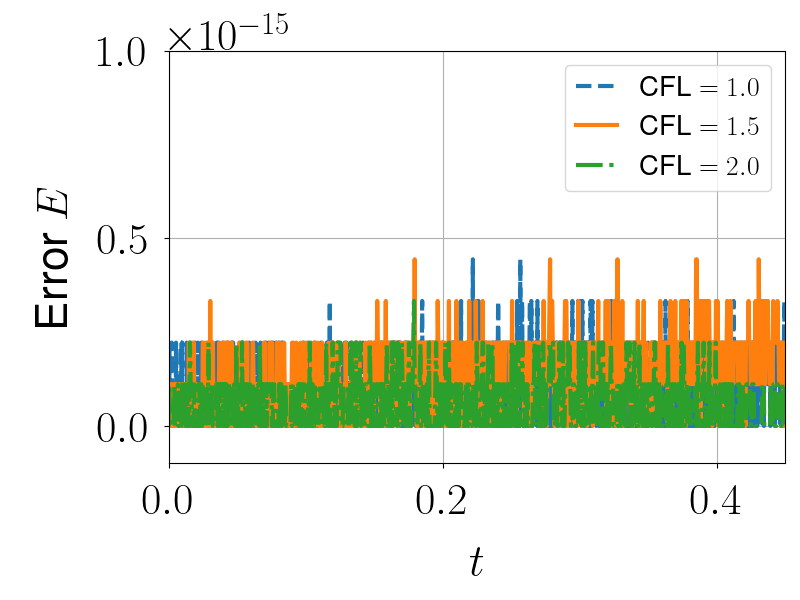}
    \end{subfigure}
    
\caption{(Mixed regime problem). Absolute error in conserved quantities as a function of time. The top and bottom rows correspond to $a_{0} = 11$ and $a_{0} = 40$, respectively.}
\label{fig:Variable Knudsen conservation data}
\end{figure}

\begin{figure}[!htbp]
    \centering
    \includegraphics[width=0.365\linewidth]{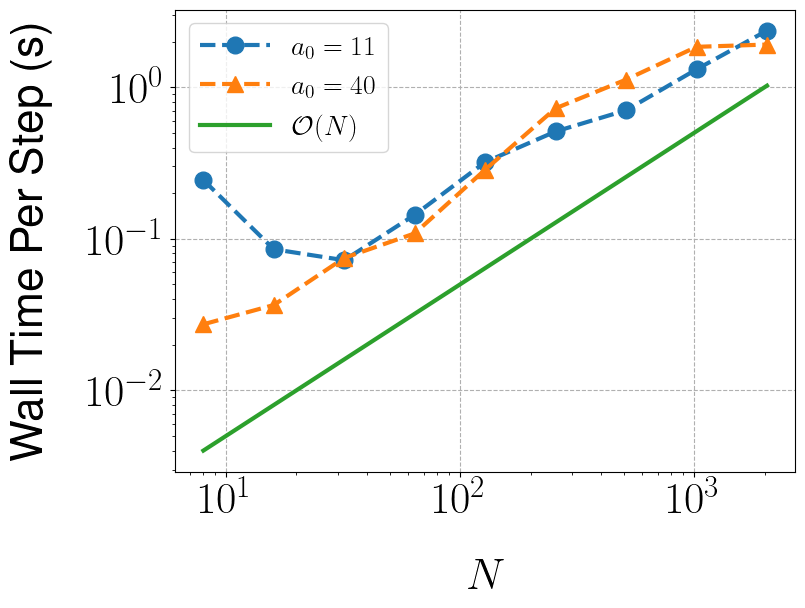}
    \caption{(Mixed regime problem). Simulation wall time per step versus the number of mesh points per dimension using a fixed CFL with $a_{0} = 11$ and $a_{0} = 40$.}
    \label{fig:Variable Knudsen scaling data}
\end{figure}

\end{exa}

\section{Conclusion}
\label{sec:conclusion}

In this work, we introduced a novel adaptive-rank method for the BGK equations that simultaneously achieves mass, momentum, and energy conservation and possesses a conditionally AP property. This method builds on our previous work by extending the greedy ACA with SVD truncation to address stiff collision operators. The greedy ACA plays a crucial role in the method and serves two key purposes: (1) It identifies the points on the mesh that are to be updated using a local SL solver and (2) enables a point-wise evaluation of the local Maxwellian in the BGK operator. The latter is particularly beneficial because local Maxwellians generally do not admit a low-rank decomposition, making this feature a distinctive advantage of our approach. In future work, we shall consider a multiplicative decomposition of the distribution function, which will expose more low-rank structure in the kinetic problem. We expect that the theorems developed in this work will be largely strengthened by such considerations, as the limiting solution has provable low-rank structure.

We also proposed a robust technique for strictly enforcing conservation laws within the kinetic solution. This technique iteratively corrects the moments of the low-rank solution using a high-order solution of the conserved macroscopic system, which is closed in a dynamic and self-consistent manner. The iterative solver for this macroscopic system is based on a highly efficient Newton-Krylov method, which allows for an entirely matrix-free implementation. Despite the fact that the provisional low-rank solution is non-conservative, it serves as a high-quality initial guess for the implicit macroscopic solver. Even for modest CFL numbers, we find that the overhead of this iteration is not significant, even when preconditioning is not applied.

To achieve high-order accuracy, we combined WENO spatial reconstructions with third-order DIRK methods that are stiffly accurate. The accuracy and efficacy of the proposed method was established using standard benchmark problems, including those with smooth and non-smooth solution structures, across a broad range of $\epsilon$. The proposed method is capable of capturing a wide range of phenomena that spans both fluid and kinetic regimes. In the final experiment, the method was applied to a challenging mixed regime problem in which the size of the Knudsen number varies by roughly six orders in magnitude. This example also highlights the advantage of the SL method, which permits a much larger time step to be used in simulations. 

In summary, the low-rank method presented in this work offers a powerful framework for solving multi-scale kinetic equations with stiff collision operators while preserving key physical properties. The flexibility of the method, combined with its efficiency and strict conservation, provides a significant advancement in the simulation of multi-scale kinetic problems. Future work will focus on extending the method to high-dimensional problems through novel low-rank formats. We also plan to explore preconditioning techniques for the macroscopic system. While simple preconditioners such as incomplete LU factorizations could be applied, they require access to the Jacobian matrix, which increases both storage and computational cost of the method. These techniques may not be feasible for high-dimensional problems; however, other strategies, such as physics-based preconditioning \cite{mousseau2000physics,knoll2001preconditioning,reisner2003efficient,hammond2005application} may be more appropriate for large scale problems in high dimensions. The BGK equation is a natural first step in the treatment of collisional problems, but it known to be accurate only when the system is ``close" to equilibrium. Future investigations will focus on more complex collision models with non-local, nonlinear effects that may better capture the full range of kinetic behaviors.

\section{Acknowledgments}

This work was partially supported Department of Energy DE-SC0023164 by the Multifaceted Mathematics Integrated Capability Centers (MMICCs) program of DOE Office of Applied Scientific Computing Research (ASCR), and Air Force Office of Scientific Research (AFOSR) FA9550-241-0254 via the Multidisciplinary University Research Initiatives (MURI) Program. We also wish to acknowledge support provided by NSF grant NSF-DMS-2111253 as well as grant FA9550-22-1-0390 from the Air Force Office of Scientific Research.

\appendix
\section{Low-Rank Numerical Integration for Macroscopic Flux}\label{app:low_rank_integration}

To efficiently evaluate the numerical integrals  $\widetilde{\int}_{\mathbb{R}^{\pm}}\cdot\,dv$, we consider the numerical integrals of an SVD solution $U\Sigma V^\top$ and an analytical Maxwellian $\mathcal{M}_{\mathcal{U}}$, which are the two basic components used to evaluate the flux \eqref{eq:BE macroscopic flux expanded}. The flux for the SVD component can be approximated using the midpoint rule, which gives
\begin{equation}
    \widetilde{\int_{\mathbb{R}^{\pm}}} f \left(\begin{array}{c}
             v\vspace{-0.2cm}  \\
             v^2\vspace{-0.2cm}\\
             \frac{1}{2}v^3
        \end{array}\right)\,dv \approx \Delta v\left(\begin{array}{l}
             U\Sigma V^\top\mathbf{v}^{\pm}\vspace{-0.2cm}  \\
             U\Sigma V^\top(\mathbf{v}^{\pm})^2\vspace{-0.2cm}\\
             \frac{1}{2}U\Sigma V^\top(\mathbf{v}^{\pm})^3
        \end{array}\right),
\end{equation}
where the wind directions are given by
\begin{equation*}
    \mathbf{v}^{+} = \max\left(\mathbf{v},0\right), \quad \mathbf{v}^{-} = \min\left(\mathbf{v},0\right).
\end{equation*}
For the Maxwellian component, we can directly use the following analytical formulas:
\begin{equation}
    \begin{split}
        &\int_{\mathbb{R}^{\pm}}\mathcal{M}_{\mathcal{U}}\,dv = \frac{\rho}{2}\left(1\pm \text{erf}\left(\frac{u}{\sqrt{2T}}\right)\right),\\
        &\int_{\mathbb{R}^{\pm}}\mathcal{M}_{\mathcal{U}}v\,dv = u\int_{\mathbb{R}^{\pm}}\mathcal{M}_{\mathcal{U}}\,dv\pm T\mathcal{M}_{\mathcal{U}}|_{v=0},\\
        &\int_{\mathbb{R}^{\pm}}\mathcal{M}_{\mathcal{U}}v^2\,dv = 2u\int_{\mathbb{R}^{\pm}}\mathcal{M}_{\mathcal{U}}v\,dv + (T-u^2)\int_{\mathbb{R}^{\pm}}\mathcal{M}_{\mathcal{U}}\,dv\mp uT\mathcal{M}_{\mathcal{U}}|_{v=0},\\
        &\int_{\mathbb{R}^{\pm}}\mathcal{M}_{\mathcal{U}}v^3\,dv = 3u\int_{\mathbb{R}^{\pm}}\mathcal{M}_{\mathcal{U}}v^2\,dv - 3u^2\int_{\mathbb{R}^{\pm}}\mathcal{M}_{\mathcal{U}}v\,dv+u^3\int_{\mathbb{R}^{\pm}}\mathcal{M}_{\mathcal{U}}\,dv+(\mp u^2\pm 2T^2 uT)\mathcal{M}_{\mathcal{U}}|_{v=0},
    \end{split}
\end{equation}
where $\rho$, $u$, and $T$ are the corresponding physical quantities of the Maxwellian $\mathcal{M}_{\mathcal{U}}$, and erf($\cdot$) is the Gauss error function.

\bibliographystyle{elsarticle-num} 
\bibliography{Reference}

\begin{thebibliography}{10}
\expandafter\ifx\csname url\endcsname\relax
  \def\url#1{\texttt{#1}}\fi
\expandafter\ifx\csname urlprefix\endcsname\relax\def\urlprefix{URL }\fi
\expandafter\ifx\csname href\endcsname\relax
  \def\href#1#2{#2} \def\path#1{#1}\fi

\bibitem{zheng2025semi}
N.~Zheng, D.~Hayes, A.~Christlieb, J.-M. Qiu, A semi-{L}agrangian adaptive-rank
  ({SLAR}) method for linear advection and nonlinear {V}lasov-{P}oisson system,
  Journal of Computational Physics 532 (2025) 113970.

\bibitem{einkemmer2025review}
L.~Einkemmer, K.~Kormann, J.~Kusch, R.~G. McClarren, J.-M. Qiu, A review of
  low-rank methods for time-dependent kinetic simulations, Journal of
  Computational Physics (2025) 114191.

\bibitem{bgk1954model}
P.~Bhatnagar, E.~Gross, M.~Krook, A model for collision processes in gases.
  {I}. small amplitude processes in charged and neutral one-component systems,
  Physical review 94~(3) (1954) 511.

\bibitem{koch2007dynamical}
O.~Koch, C.~Lubich, Dynamical low-rank approximation, SIAM Journal on Matrix
  Analysis and Applications 29~(2) (2007) 434--454.

\bibitem{einkemmer2018low}
L.~Einkemmer, C.~Lubich, A low-rank projector-splitting integrator for the
  {Vlasov}--{Poisson} equation, SIAM Journal on Scientific Computing 40~(5)
  (2018) B1330--B1360.

\bibitem{einkemmer2019low}
L.~Einkemmer, A low-rank algorithm for weakly compressible flow, SIAM Journal
  on Scientific Computing 41~(5) (2019) A2795--A2814.

\bibitem{einkemmer2021efficient}
L.~Einkemmer, J.~Hu, L.~Ying, An efficient dynamical low-rank algorithm for the
  {B}oltzmann-{BGK} equation close to the compressible viscous flow regime,
  SIAM Journal on Scientific Computing 43~(5) (2021) B1057--B1080.

\bibitem{ceruti2022unconventional}
G.~Ceruti, C.~Lubich, An unconventional robust integrator for dynamical
  low-rank approximation, BIT Numerical Mathematics 62~(1) (2022) 23--44.

\bibitem{baumann2024stable}
L.~Baumann, L.~Einkemmer, C.~Klingenberg, J.~Kusch, A stable multiplicative
  dynamical low-rank discretization for the linear {B}oltzmann-{BGK} equation,
  arXiv preprint arXiv:2411.06844 (2024).

\bibitem{kormann2015semi}
K.~Kormann, A semi-{Lagrangian} {Vlasov} solver in tensor train format, SIAM
  Journal on Scientific Computing 37~(4) (2015) B613--B632.

\bibitem{dektor2021rank}
A.~Dektor, A.~Rodgers, D.~Venturi, Rank-adaptive tensor methods for
  high-dimensional nonlinear {PDEs}, Journal of Scientific Computing 88~(2)
  (2021) 36.

\bibitem{GuoVlasovFlowMap2022}
W.~Guo, J.-M. Qiu, {A low rank tensor representation of linear transport and
  nonlinear Vlasov solutions and their associated flow maps}, Journal of
  Computational Physics 458 (2022) 111089.
\newblock \href {https://doi.org/111089} {\path{doi:111089}}.

\bibitem{guo2024local}
W.~Guo, J.~F. Ema, J.-M. Qiu, A {Local} {Macroscopic} {Conservative} {(LoMaC)}
  low rank tensor method with the discontinuous {Galerkin} method for the
  {Vlasov} dynamics, Communications on Applied Mathematics and Computation
  6~(1) (2024) 550--575.

\bibitem{sands2025high}
W.~A. Sands, W.~Guo, J.-M. Qiu, T.~Xiong, High-order adaptive rank integrators
  for multi-scale linear kinetic transport equations in the hierarchical
  {T}ucker format, {SIAM} {J}ournal on {S}cientific {C}omputing (to appear)
  (2025).

\bibitem{shi2024distributed}
T.~Shi, D.~Hayes, J.-M. Qiu, Distributed memory parallel adaptive tensor-train
  cross approximation, arXiv preprint arXiv:2407.11290 (2024).

\bibitem{dektor2025nonlinear}
A.~Dektor, Collocation methods for nonlinear differential equations on low-rank
  manifolds, Linear Algebra and its Applications 705 (2025) 143--184.

\bibitem{dektor2024interpolatoryBGK}
A.~Dektor, L.~Einkemmer, Interpolatory dynamical low-rank approximation for the
  3+ 3d {B}oltzmann-{BGK} equation, arXiv preprint arXiv:2411.15990 (2024).

\bibitem{guo2024conservative}
W.~Guo, J.-M. Qiu, A conservative low rank tensor method for the {Vlasov}
  dynamics, SIAM Journal on Scientific Computing 46~(1) (2024) A232--A263.

\bibitem{li2023high}
L.~Li, J.~Qiu, G.~Russo, A high-order {semi-Lagrangian finite difference method
  for nonlinear Vlasov and {BGK} models}, Communications on Applied Mathematics
  and Computation 5~(1) (2023) 170--198.

\bibitem{ding2023accuracy}
M.~Ding, J.-M. Qiu, R.~Shu, Accuracy and stability analysis of the
  semi-lagrangian method for stiff hyperbolic relaxation systems and kinetic
  {BGK} model, SIAM Multiscale Modeling \& Simulation 21~(1) (2023) 143--167.

\bibitem{taitano2014moment}
W.~T. Taitano, D.~A. Knoll, L.~Chac{\'o}n, J.~M. Reisner, A.~K. Prinja,
  Moment-based acceleration for neutral gas kinetics with {BGK} collision
  operator, Journal of Computational and Theoretical Transport 43~(1-7) (2014)
  83--108.

\bibitem{taitano2015charge}
W.~T. Taitano, D.~A. Knoll, L.~Chac{\'o}n, Charge-and-energy conserving
  moment-based accelerator for a multi-species
  {V}lasov--{F}okker--{P}lanck--{A}mp{\`e}re system, part {II}: {C}ollisional
  aspects, Journal of Computational Physics 284 (2015) 737--757.

\bibitem{park2019multigroup}
H.~Park, L.~Chac{\'o}n, A.~Matsekh, G.~Chen, A multigroup moment-accelerated
  deterministic particle solver for 1-{D} time-dependent thermal radiative
  transfer problems, Journal of Computational Physics 388 (2019) 416--438.

\bibitem{hammer2019multi}
H.~Hammer, H.~Park, L.~Chac{\'o}n, A multi-dimensional, moment-accelerated
  deterministic particle method for time-dependent, multi-frequency thermal
  radiative transfer problems, Journal of Computational Physics 386 (2019)
  653--674.

\bibitem{chapman1990mathematical}
S.~Chapman, T.~G. Cowling, The mathematical theory of non-uniform gases: an
  account of the kinetic theory of viscosity, thermal conduction and diffusion
  in gases, Cambridge university press, 1990.

\bibitem{civril2007finding}
A.~Civril, M.~Magdon-Ismail, Finding maximum volume sub-matrices of a matrix,
  RPI Comp Sci Dept TR (2007) 07--08.

\bibitem{1994Kinetic}
J.~C. Mandal, S.~M. Deshpande, Kinetic flux vector splitting for {Euler}
  equations, Computers \& Fluids 23~(2) (1994) 447--478.

\bibitem{1995Gas}
K.~Xu, L.~Martinelli, A.~Jameson, Gas-kinetic finite volume methods,
  flux-vector splitting, and artificial diffusion, Journal of Computational
  Physics 120~(1) (1995) 48--65.

\bibitem{JIANG1996202}
G.-S. Jiang, C.-W. Shu, Efficient implementation of weighted {ENO} schemes,
  Journal of Computational Physics 126~(1) (1996) 202--228.

\bibitem{knoll2004jacobian}
D.~A. Knoll, D.~E. Keyes, Jacobian-free {N}ewton--{K}rylov methods: a survey of
  approaches and applications, Journal of Computational Physics 193~(2) (2004)
  357--397.

\bibitem{knoll2005jacobian}
D.~A. Knoll, V.~Mousseau, L.~Chac{\'o}n, J.~Reisner, Jacobian-free
  {N}ewton-{K}rylov methods for the accurate time integration of stiff wave
  systems, Journal of Scientific Computing 25 (2005) 213--230.

\bibitem{alexander1977diagonally}
R.~Alexander, Diagonally implicit {Runge--Kutta} methods for stiff {ODE’s},
  SIAM Journal on Numerical Analysis 14~(6) (1977) 1006--1021.

\bibitem{kennedy2016diagonally}
C.~A. Kennedy, M.~H. Carpenter, Diagonally implicit {R}unge-{K}utta methods for
  ordinary differential equations. a review, Tech. rep., NASA Langley Research
  Center (2016).

\bibitem{pieraccini2007implicit}
S.~Pieraccini, G.~Puppo, Implicit--explicit schemes for {BGK} kinetic
  equations, Journal of Scientific Computing 32 (2007) 1--28.

\bibitem{xiong2015high}
T.~Xiong, J.~Jang, F.~Li, J.-M. Qiu, High order asymptotic preserving nodal
  discontinuous galerkin {IMEX} schemes for the {BGK} equation, Journal of
  Computational Physics 284 (2015) 70--94.

\bibitem{mousseau2000physics}
V.~Mousseau, D.~Knoll, W.~Rider, Physics-based preconditioning and the
  {N}ewton--{K}rylov method for non-equilibrium radiation diffusion, Journal of
  computational physics 160~(2) (2000) 743--765.

\bibitem{knoll2001preconditioning}
D.~A. Knoll, W.~Vanderheyden, V.~Mousseau, D.~B. Kothe, On preconditioning
  {N}ewton--{K}rylov methods in solidifying flow applications, SIAM Journal on
  Scientific Computing 23~(2) (2001) 381--397.

\bibitem{reisner2003efficient}
J.~Reisner, A.~Wyszogrodzki, V.~Mousseau, D.~Knoll, An efficient physics-based
  preconditioner for the fully implicit solution of small-scale thermally
  driven atmospheric flows, Journal of Computational Physics 189~(1) (2003)
  30--44.

\bibitem{hammond2005application}
G.~E. Hammond, A.~J. Valocchi, P.~C. Lichtner, Application of {J}acobian-free
  {N}ewton--{K}rylov with physics-based preconditioning to biogeochemical
  transport, Advances in Water Resources 28~(4) (2005) 359--376.

\end{thebibliography}

\end{document}